\documentclass[reprint, aps, twocolumn,groupedaddress, showkeys]{revtex4-2}

\usepackage[dvipsnames]{xcolor}
\usepackage{soul}
\usepackage{latexsym,amsmath,amssymb,graphics,amscd, empheq}
\usepackage{graphicx}
\usepackage{float}
\usepackage{subcaption}
\usepackage{color}
\usepackage{stackengine}
\usepackage[breaklinks,colorlinks,
  citecolor=blue,
  urlcolor=blue]{hyperref}

\usepackage{enumitem}

\newtheorem{theorem}{Theorem}[section]

\newtheorem{lemma}[theorem]{Lemma}
\newtheorem{remark}[theorem]{Remark}
\newtheorem{proposition}[theorem]{Proposition}

\providecommand{\norm}[1]{\lVert#1\rVert}

\newcommand{\vertiii}[1]{{\left\vert\kern-0.25ex\left\vert\kern-0.25ex\left\vert #1 
    \right\vert\kern-0.25ex\right\vert\kern-0.25ex\right\vert}}
\newsavebox{\savepar}

\makeatletter

\usepackage{scalerel,stackengine}
\stackMath
\newcommand\reallywidehat[1]{%
\savestack{\tmpbox}{\stretchto{%
  \scaleto{%
    \scalerel*[\widthof{\ensuremath{#1}}]{\kern-.6pt\bigwedge\kern-.6pt}%
    {\rule[-\textheight/2]{1ex}{\textheight}}
  }{\textheight}%
}{0.5ex}}%
\stackon[1pt]{#1}{\tmpbox}%
}

\begin{document}


\title{{Chaos on compact manifolds:
Differentiable synchronizations beyond the Takens theorem}}


\author{Lyudmila Grigoryeva}
\email[Lyudmila.Grigoryeva@uni-konstanz.de]{}
\affiliation{Department of Mathematics and Statistics. Universit\"at Konstanz. Box 146. D-78457 Konstanz. Germany.}

\author{Allen Hart}
\email[A.Hart@bath.ac.uk]{}
\affiliation{Department of Mathematical Sciences, University of Bath, Bath BA2 7AY, UK.}

\author{Juan-Pablo Ortega}
\email[Juan-Pablo.Ortega@ntu.edu.sg]{}
\affiliation{Division of Mathematical Sciences. Nanyang Technological University.
Singapore.} 



\begin{abstract}
This paper shows that a large class of fading memory state-space systems driven by discrete-time observations of dynamical systems  defined on compact manifolds always yields continuously differentiable synchronizations.  This general result provides a powerful tool for the representation, reconstruction, and forecasting of chaotic attractors. It also improves previous statements in the literature for differentiable generalized synchronizations, whose existence was so far  guaranteed for a restricted family of systems and was detected using H\"older exponent-based criteria.
\end{abstract}

\keywords{dynamical systems, synchronization, chaos, attractor, Takens embedding, echo state property, fading memory property, asymptotic stability, echo state network.}

\maketitle

\section{\label{Introduction}Introduction}

Synchronization phenomena between chaotic systems have deserved much attention for decades (see \cite{pecora:synch, ott2002chaos, boccaletti:reports:2002, eroglu2017synchronisation} for self-contained presentations and many references). This topic has intrinsic theoretical interest {and} a variety of applications ranging from encryption \cite{alvarez2005breaking, moskalenko2010generalized} and communication schemes \cite{jiang2004chaotic} to the analysis of neurological disorders \cite{stam2002generalized, bartolomei2006disturbed}.  {\it Generalized synchronizations} (GS) were introduced in \cite{rulkov1995generalized} to {describe}  a system driven by the observations of a chaotic dynamical system that (asymptotically) yields a master-slave configuration. The master is the full dynamical system (not its observations) and the slave is the driven system. Any map that implements this configuration is called a GS.

Given a system driven by the observations of an invertible dynamical system,  the main theorem in \cite{kocarev1995general} shows that the asymptotic stability of the system response is a sufficient condition for the existence of a GS. Nevertheless, it was quickly noticed in \cite{pyragas:1996, hunt:ott:1997} that the GS guaranteed by this theorem may have poor regularity properties, rendering it useless  as an attractor representation and reconstruction tool. {This fact motivated the distinction in \cite{pyragas:1996} between {\it strong} and {\it weak} GS depending on whether the synchronization map is differentiable or not. Differentiability was also identified in \cite{hunt:ott:1997} as a crucial property that determines how useful a GS may be in important tasks such as the estimation of attractor dimensions or Lyapunov exponents. Additionally, a differentiability criterion was introduced in that paper based on the H\"older exponent of the response, mostly for invertible driving dynamical systems whose attractors have only one negative Lyapunov exponent.}

{In this paper we replace the analytical approach in \cite{hunt:ott:1997} with a geometrical one that, combined with arguments coming from functional analysis, allows us to prove another differentiability criterion. Our main statement (Theorem \ref{main theorem of letter}) shows that if  the invertible dynamical system} {\it has a compact manifold as phase space and its observations drive a fading memory system} with locally state-contracting map, then {\it  the generalized synchronizations introduced in \cite{kocarev1995general} always exist and are continuously differentiable}. We carefully define all these terms later on in the paper.


Recurrent neural networks and reservoir computing \cite{lukosevicius, tanaka:review} have recently enjoyed remarkable success  learning \cite{Jaeger04, pathak:chaos, Pathak:PRL, Ott2018} and classifying \cite{carroll2018using} chaotic attractors of complex nonlinear infinite dimensional dynamical systems, and detecting GS phenomena \cite{ibanez2018detection, Weng:2019, Lymburn:2019}. This strongly suggested that these machine learning paradigms have Takens embedding-type properties   \cite{takensembedding, huke:2006}. This fact has been rigorously established in \cite{hart:ESNs} where the so called {\it Echo State Networks} (ESNs) \cite{Matthews:thesis, Matthews1993, Jaeger04, RC7, RC8} driven by one-dimensional observations of a given dynamical system on a compact manifold have been shown to produce dynamics that are topologically conjugate to that of the original system. This result is actually proved by, among other things, showing that a natural map that arises in that context (called Echo State Map) is a differentiable generalized synchronization.

Theorem \ref{main theorem of letter} in this paper is an {extensive generalization of that differentiability statement, that is shown to be valid not only for ESNs}, but for {\it any fading memory system generated by a locally state-contracting system}.

{In many of the most recent studies (as in \cite{Jaeger04, pathak:chaos, Pathak:PRL, Ott2018, carroll2018using}), learning and classifying chaotic attractors is tackled by constructing state-space systems to which one naturally associates autonomous systems that serve as proxies for the dynamics that needs to be learnt. An important implementation problem that arises using this approach is the so-called {\it curse of dimensionality} which, in general, can be elegantly avoided assuming higher order regularity \cite{Poggio2017}. More specifically,  the availability of an injective GS for a dynamical system implies the \emph{supervised learnability} (explained in the next section)} {of the dynamics from} {its observations. The learning is done by approximating a readout function that is defined in a potentially high-dimensional state space. We shall see that the regularity of this readout is determined by the smoothness of the dynamical system in question and of the inverse of the GS (available by the injectivity hypothesis). Hence, if this readout is learnt using standard paradigms like (deep) neural networks or splines, then classical works (see for instance \cite{Poggio2017, Mhaskar1996}) show that the approximation rate depends on the ratio between the state space dimension (again, potentially very high) and  the smoothness of the readout. The differentiability of the GS is therefore of crucial practical importance.}

A more recent framework in which GS differentiability arises as a key feature is that of the {\it transfer learning} of dynamical systems, as characterized in \cite{inubushi2020transfer}. This technique translates to the dynamical systems context a concept that has already been successfully used in machine learning \cite{weiss2016survey} and that spells out how to adequately learn the dynamics of a system with training data generated by another one that is {\it close} in systems space. When the different systems to be learnt are labeled by a parameters space, a necessary condition for this technique to function is the structural stability of the data generating systems with respect to those parameters, as well as the differentiability of the readouts that are trained to implement the learning. As we explain in the next section, those readouts can be encoded in terms of a GS whose differentiability ensures that of the corresponding readout and hence makes transfer learning possible.

{In view of this discussion, this paper identifies a rich class of systems  for which continuously differentiable GS are available} {for invertible dynamical systems on compact manifolds.} {Additionally, when the GS happens to be injective, the supervised learning of the attractor of a given dynamical system} {based on finite-length observations is feasible in practice. These results are valuable for the representation, reconstruction, and forecasting of chaotic attractors and explain the excellent performance of }{reservoir computing and of ESNs in particular evidenced in the above references.}

\section{\label{Definitions and preliminary discussion}Definitions and preliminary discussion}

\noindent {\bf The dynamical system.}
All along this paper we consider an invertible and discrete-time dynamical system determined by a map  $\phi$ that belongs either to the set  ${\rm Hom}(M)$ of homeomorphisms (continuous invertible maps with continuous inverse) of a compact topological space $M$ or to the set of diffeomorphisms  ${\rm Diff} ^1(M) $ of a finite-dimensional compact manifold $M $. Later on, in the differentiable case, we need to ensure that $M$ can be endowed with a Riemannian metric and that is why we additionally assume that $M$ is connected, Hausdorff, and second-countable (see \cite[Proposition 2.10]{do:carmo:1993}). The $d$-dimensional observations of the dynamical system are realized by maps $\omega $ that belong either to $C ^0(M, \mathbb{R}^d)$ or to $C ^1(M, \mathbb{R}^d) $ depending on the nature of $\phi  $ (${\rm Hom}(M)$ or  ${\rm Diff} ^1(M) $, respectively). When $M$ is a manifold, the symbol $TM $ denotes the tangent bundle of  $M$, $T \phi:TM \longrightarrow TM $ the tangent map of $\phi $,  and $D \omega: TM \longrightarrow \mathbb{R}^d   $  the differential of the observation map $\omega $.

\medskip

\noindent {\bf Sequences.}
$\Bbb Z $ denotes the integers and the symbol $\mathbb{Z}_{-}:= \left\{\ldots, -2,-1,0\right\} $   stands for the non-positive integers. Given $D_d \subset \mathbb{R} ^d $, we denote by $D_d ^{\Bbb Z}$ (respectively $D_d^{\mathbb{Z}_{-} } $) the set of $D_d$-valued two-sided  infinite (respectively semi-infinite) sequences. The symbol $(\ell^{\infty}(\mathbb{R}^d), \left\|\cdot \right\|_{\infty}) $ (respectively $(\ell^{\infty}_-(\mathbb{R}^d), \left\|\cdot \right\|_{\infty}) $) denotes the Banach space of $\mathbb{R} ^d $-valued two-sided  infinite (respectively semi-infinite) sequences that have a finite supremum. For any $t \in \Bbb Z$ we define the {\it projection}  $p _t:\ell^{\infty}(\mathbb{R}^d) \longrightarrow \mathbb{R}^d $ such that $p _t({\bf z}):= {\bf z} _t  $ and the {\it time delay} operator $T _t:\ell^{\infty}(\mathbb{R}^d) \longrightarrow \ell^{\infty}(\mathbb{R}^d) $ given by $T_t( {\bf z}) _\tau:= {\bf z}_{\tau- t}$, $\tau \in \Bbb Z  $. It is easy to see that both $p _t $ and $T _t $ are bounded linear operators and that for any $t _1, t _2  \in \mathbb{Z}$  we have 
\begin{equation}
\label{relation p and T}
p_{t _1+ t _2}= p_{t _1} \circ T_{-t _2}=p_{t _2} \circ T_{-t _1}.
\end{equation}
These definitions can be easily adapted to time delays and projections defined in $\ell^{\infty}_-(\mathbb{R}^d) $ (see \cite[Section 2.3]{RC9}).

\medskip

\noindent {\bf The delay map.}
Define the $(\phi, \omega)$-{\it delay map} $S _{(\phi, \omega)}:M \longrightarrow \ell^{\infty}(\mathbb{R}^d) $ as $S _{(\phi, \omega)}(m):=\left\{\omega(\phi ^t (m))\right\}_{t \in \Bbb Z} $ and $S _{(\phi, \omega)}^-:M \longrightarrow \ell^{\infty}_-(\mathbb{R}^d) $ by $S _{(\phi, \omega)}^-(m):=\left\{\omega(\phi ^t (m))\right\}_{t \in \Bbb Z_-} $. Note that the Abelian group $(\Bbb Z, +) $ acts both on the phase space $M$ via  $\phi $ and on the  space $\ell^{\infty}(\mathbb{R}^d) $ via the time delay operators. The map $S _{(\phi, \omega)} $ is equivariant with respect to those two actions, that is,
\begin{equation}
\label{somega equivariant}
T_{-t}(S _{(\phi, \omega)}(m))= S _{(\phi, \omega)}(\phi ^t(m)),
\end{equation}
for all $t \in \Bbb Z$ and $m \in M$.

\medskip

\noindent {\bf State maps, the echo state property (ESP), and the fading memory property (FMP).}
Let  $F: \mathbb{R} ^N\times \mathbb{R} ^d \longrightarrow   \mathbb{R} ^N $ be a state map with states in $\mathbb{R}^N  $, $ N \in \mathbb{N} $. Whenever there exist subsets $D_N \subset \mathbb{R} ^N  $  and $D _d\subset \mathbb{R} ^d $ such that $F(D _N\times D _d)\subset D _N  $, we shall  denote the restricted state map $F:D _N\times D _d\longrightarrow D _N  $ with the same symbol. Given a state map $F:D _N\times D _d\longrightarrow D _N  $ and two subsets $V _d\subset D_d^{\mathbb{Z}} $ and $V _N\subset D_N^{\mathbb{Z}} $, we say that the state map determined by $F$ has the $(V _d, V _N) $-{\it echo state property} (ESP) (see \cite{jaeger2001, Manjunath:Jaeger, manjunath:prsl} for in-depth descriptions of this property) when for any ${\bf z} \in V _d\subset D_d^{\mathbb{Z}} $ there exists a unique $\mathbf{x}\in V _N\subset D_N^{\mathbb{Z}} $ such that 
\begin{equation}
\label{state equations for ESP}
\mathbf{x} _t=F(\mathbf{x} _{t-1}, {\bf z} _t), \quad \mbox{for all} \quad t \in \mathbb{Z}.
\end{equation}
 
When the inputs in the system determined by $F:D _N\times D _d\longrightarrow D _N  $ are the $\omega $-observations of the dynamical system $\phi$ (we assume that $\omega(\phi(M))\subset D _d$) we then talk about the $(\phi, \omega) $-{\it echo state property}. This is a particular case of the general definition that we just introduced as  the $(\phi, \omega) $-ESP is the same as the $(S _{(\phi, \omega)}(M), D _N^{\mathbb{Z}}) $-ESP in the sense above.

Given a system like \eqref{state equations for ESP} that has the $(V _d, V _N) $-ESP we can naturally define a causal and time-invariant  filter $U ^F: V _d \longrightarrow V _N  $ by assigning to each ${\bf z} \in V _d $ the unique $\mathbf{x} \in V _N $ that satisfies \eqref{state equations for ESP}. We recall that the dynamics of causal and time-invariant filters  is fully determined by their restriction (in the domain and codomain) to semi-infinite sequences labeled by  $ \mathbb{Z}_{-}$. We will use exchangeably  the filters  operating on two-sided infinite and semi-infinite sequences. 

Let $U: \ell^{\infty}_-(\mathbb{R}^d) \longrightarrow \ell^{\infty}_-(\mathbb{R}^N) $ be a causal and time-invariant filter. Consider a {\it  weighting sequence} $w$,
that is, a strictly decreasing sequence with zero limit  $w : \mathbb{N} \longrightarrow (0,1] $ such that $w _0=1 $, and define $w${\it -norm} by
\begin{equation*}
\| {\bf z} \| _{w}:= \sup_{t \in \Bbb Z_-}\{\| {\bf z}_t\| w_{-t}\},\quad \mbox{for any} \quad {\bf z}\in ({\Bbb R}^d)^{\mathbb{Z}_{-}}.
\end{equation*}
We say that $U$ has the {\it  fading memory property} (FMP)
with respect to the sequence $w$ if $U: (\ell^{\infty}_-(\mathbb{R}^d), \norm{.}_w) \longrightarrow ( \ell^{\infty}_-(\mathbb{R}^N) , \norm{.}_w)$ is continuous, that is for any ${\bf z} \in \ell^{\infty}_-(\mathbb{R}^d)$ and any $\epsilon > 0$ there exists a $\delta(\epsilon, {\bf z}) > 0$
such that for any $\bar{\bf z} \in \ell^{\infty}_-(\mathbb{R}^d)$ it holds that
\begin{equation}
\label{definition FMP rc18}
  \norm{\bar{\bf z} - {\bf z}}_w < \delta(\epsilon, {\bf z}) \implies \norm{U(\bar{\bf z}) - U({\bf z})}_w < \epsilon.
\end{equation}
We recall that for any compact set $D _d  \subset \mathbb{R} ^d $ (see \cite[Corollary 2.7]{RC6}) the metric induced by any weighted norm $\left\|\cdot \right\|_w $  on $D_d^{\mathbb{Z}_{-}}\subset \ell_{-}^{\infty}(\mathbb{R}^d)$ generates the product topology on it. The FMP is a dynamic continuity feature that, roughly speaking, makes the influence of the inputs on the outputs of the system determined by $F$ less important  as they become more distant in the past. 

We shall use the following notation: in Euclidean spaces $\mathbb{R}^d $ and for any  $L>0 $, the symbol $B _L \subset \mathbb{R}^d $ (respectively $B _L (\mathbf{v}) \subset \mathbb{R}^d $, with $\mathbf{v} \in \mathbb{R}^d $)  denotes the  $2$-ball centered at zero (respectively at $\mathbf{v} \in {\Bbb R}^d $). 

The next result shows that {\it locally contracting state maps produce systems that have automatically the ESP and the FMP}, provided that the {\it inputs take values in a compact set}.

\begin{proposition}
\label{contracting implies omega ESP}
Let  $F: \mathbb{R} ^N\times \mathbb{R} ^d \longrightarrow   \mathbb{R} ^N $ be a continuous state map,   let $D _N \subset \mathbb{R} ^N  $ be a closed subset, and let $D _d \subset \mathbb{R} ^d  $ be a compact subset. Suppose that $F(D _N \times D_d )\subset D _N $ and that the restriction (denoted with the same symbol) $F:D _N \times D_d \longrightarrow  D _N $ is a contraction in the state variable, that is, there exists a constant $0<c<1 $ such that for all $\mathbf{x} _1, \mathbf{x} _2 \in D _N$ and ${\bf z} \in D_d  $
\begin{equation}
\label{global contractivity hypothesis}
\|F({\bf x}_1,{\bf z})- F({\bf x}_2,{\bf z}) \| \leq c \|{\bf x}_1-{\bf x}_2\|. 
\end{equation}
Then:
\begin{description}
\item [(i)] There exists a compact subset $W \subset D _N $ such that $F(W\times D_d) \subset W$ and the system determined by $F: W\times D_d \longrightarrow   W$ with inputs in $D _d^{\mathbb{Z}_{-}}\subset \ell^{\infty}_-(\mathbb{R}^d) $ has the $(D _d^{\mathbb{Z}_{-}}, W ^{\mathbb{Z}_{-}})$-ESP. This means that for any ${\bf z} \in D _d^{\mathbb{Z}_{-}}$ there exists a unique $\mathbf{x}\in W ^{\mathbb{Z}_{-}}$ such that 
\begin{equation}
\label{state system eq}
\mathbf{x} _t=F( \mathbf{x} _{t-1}, {\bf z} _t), \quad \mbox{for all} \quad t \in \Bbb Z_-.
\end{equation}
If $D_N$ is convex then the compact subset $W$ can be chosen  to be also convex.
\item [(ii)] The recursions \eqref{state system eq} determine uniquely a state filter $U^F:D _d^{\mathbb{Z}_{-}} \longrightarrow W ^{\mathbb{Z}_{-}}$ that satisfies 
\begin{equation}
\label{state system eq filter}
U^F({\bf z}) _t=F\left(U^F({\bf z}) _{t-1}, {\bf z} _t\right), \quad \mbox{for all} \quad t \in \Bbb Z_-.
\end{equation}
The filter $U^F $ is continuous when $D _d^{\mathbb{Z}_{-}}$  and  $W ^{\mathbb{Z}_{-}}$ are endowed with the product topology. Moreover, $U^F $ has the fading memory property with respect to any weighting sequence.
\end{description}
\end{proposition}

\medskip

\noindent {\bf Generalized synchronizations.} Consider now a continuously differentiable state map $F: \mathbb{R} ^N\times \mathbb{R} ^d \longrightarrow   \mathbb{R} ^N $ with states in $\mathbb{R}^N  $, $ N \in \mathbb{N} $, as well as the  drive-response system associated to the $\omega $-observations of the dynamical system $\phi$ and determined by the recursions:
\begin{equation}
\label{drive-response system}
\mathbf{x} _t=F\left(\mathbf{x} _{t-1}, S _{(\phi, \omega)}(m)_t\right), \quad \mbox{$t \in \Bbb Z,\, m \in M.$}
\end{equation}
We say that a {\it generalized synchronization} (GS) occurs in this configuration when there exists a map $f:M \longrightarrow \mathbb{R}^N $ such that for any $\mathbf{x} _t  $, $t \in \Bbb Z $, as in \eqref{drive-response system} it holds that
\begin{equation}
\label{generalized synchronization condition}
 \mathbf{x} _t = f (\phi ^t(m)),
\end{equation}
that is, the time evolution of the dynamical system in phase space (not just its observations) drives the response in \eqref{drive-response system}. All these concepts can be easily extended to the more general situation in which $M$ is just a compact topological space, $\phi \in {\rm Hom} (M) $ is a homeomorphism (that is, $\phi $ is continuous, invertible, and the inverse is also continuous), and $\omega$ and $F$ are just continuous.

When a GS is invertible, the readout $h:= \omega \circ \phi \circ f ^{-1}: f(M) \subset  \mathbb{R}^N \longrightarrow \mathbb{R} ^d $ of the states $\mathbf{x}_t $ determined by \eqref{generalized synchronization condition} fully characterize the dynamics of the $\omega $-observations $\left\{\omega(\phi ^t (m))\right\}_{t \in \Bbb Z} $ of $\phi$ because $ h(\mathbf{x} _t)= \omega\circ  \phi(f ^{-1}(\mathbf{x} _t))= \omega(\phi ^{t+1} (m))$. This observation implies that this dynamics can be captured via supervised learning if the function $h:= \omega \circ \phi \circ f ^{-1} $ is sufficiently regular to be efficiently approximated by a universal family (for instance neural networks or polynomials). This rationale has been profusely exploited in forecasting applications in the context of reservoir computing (see  \cite{Jaeger04, pathak:chaos, Pathak:PRL, Ott2018, carroll2018using} as well as \cite{lu:bassett:2020, Verzelli2020b}). We emphasize that in the reservoir computing context the state map $F$ is quasi-randomly generated and the readout $h$ is constrained to be linear. We are not aware of the existence of a dynamic representation theorem of the type that we just evoked in the presence of invertible synchronizations and we conjecture that such a statement can only be formulated up to topological conjugacies, in the spirit of the results in \cite{hart:ESNs, allen:tikhonov}. 

We emphasize that the definition \eqref{generalized synchronization condition} presupposes that $F$ has the $(\phi, \omega) $-echo state property and that $F$ hence determines a unique causal and time-invariant filter $U ^F: S _{(\phi, \omega)}(M) \longrightarrow (\mathbb{R} ^N)^{\Bbb Z} $ that associates to each orbit $S _{(\phi, \omega)}(m) $ the unique solution sequence $\mathbf{x} \in (\mathbb{R} ^N)^{\Bbb Z} $ of \eqref{drive-response system}.

In the following lemma we show that
the map $f_{(\phi, \omega,F)} :M \longrightarrow \mathbb{R}^N $ defined by 
\begin{equation}
\label{definition our gs}
f_{(\phi, \omega,F)} (m):=p _0 \left(U ^F(S _{(\phi, \omega)} (m))\right), 
\end{equation}
with $p _0:(\mathbb{R} ^N)^{\Bbb Z} \rightarrow \mathbb{R} ^N $ the projection onto the zero entry of the sequence, is a GS between the original dynamical system $\phi $ and the response of the system $F$ driven by its $\omega$-observations. 

\begin{lemma}
\label{our f is a GS}
Let $\phi \in {\rm Hom} (M) $ be an invertible dynamical system on a compact topological space $M $, $\omega \in C ^0(M, \mathbb{R}^d) $ an observation map, and $F: D _N\times D _d \longrightarrow   D_N $ a continuous state map, with $D _N\subset \mathbb{R} ^N $ and $D_d\subset \mathbb{R}^d  $. If the system determined by $F$ and driven by the $\omega $-observations of $\phi $ has the $(\phi, \omega) $-ESP, then the map $f_{(\phi, \omega,F)} :M \longrightarrow D_N $ defined by $f_{(\phi, \omega,F)} (m):=p _0 \left(U ^F(S _{(\phi, \omega)} (m))\right) $ is a generalized synchronization, that is, it satisfies the defining relation \eqref{generalized synchronization condition} and, more generally,
\begin{equation}
\label{defining claim for GS}
U^F(S _{(\phi, \omega)}(m)) _t=f_{(\phi, \omega,F)} \left(\phi ^t (m)\right), 
\end{equation}
for any $t \in \Bbb Z,\, m \in M$.
\end{lemma}

\begin{lemma}
\label{lemma SSM recursion}
In the same setup as Lemma \ref{our f is a GS}, the state synchronization map $f_{(\phi, \omega,F)} $ satisfies the identity:
\begin{equation}
\label{SSM recursion for lemma}
f_{(\phi, \omega,F)} (m)=F \left(f_{(\phi, \omega,F)} (\phi ^{-1}(m)), \omega (m)\right), 
\end{equation}
for all $m \in  M $. 
\end{lemma}

In the next section, the main result of the paper, Theorem \ref{main theorem of letter}, characterizes a large family of state maps $F$ for which the GS $f_{(\phi, \omega,F)}  $ {\it is continuously differentiable} and, additionally, {\it it is the unique generalized synchronization that satisfies the recursion} \eqref{SSM recursion for lemma}. More specifically, the crucial condition on the state map $F$ that makes GS of the type \eqref{definition our gs} available is its {\it local contractivity in the state variable}. This means that there exists a constant $0<c<1 $ and a closed convex set $V \subset \mathbb{R}^N  $ such that $F(V \times \omega(M))\subset V $ and, additionally, for all $\mathbf{x} _1, \mathbf{x} _2 \in V$ and ${\bf z} \in \omega (M)$,
\begin{equation}
\label{eq:Fcontractive paper} 
\|F({\bf x}_1,{\bf z})- F({\bf x}_2,{\bf z}) \| \leq c \|{\bf x}_1-{\bf x}_2\|. 
\end{equation}
This condition ensures (see Propositions \ref{contracting implies omega ESP} and \ref{continuous synchronization}) that the ESP always holds. Additionally, we recall that the local contractivity implies the  fading memory property (see \eqref{definition FMP rc18} and part {\bf (ii)} in Proposition \ref{contracting implies omega ESP}). The FMP implies the so called {\it unique steady-state property}, also referred to as {\it input-forgetting property} \cite{jaeger2001, RC9}. All these properties are relevant in our context as they all coincide with the {\it asymptotic stability} that has been identified in the foundational paper \cite{kocarev1996generalized} as the characterizing property for the existence of generalized synchronizations.

\section{The main result}

The next theorem, proved in detail in the appendices  (see Proposition \ref{continuous synchronization}, Theorem \ref{differentiable generalized synchronizations theorem}, and Remark \ref{global vs lokal 1}), shows that for invertible dynamical systems $\phi$ defined on compact topological spaces $M$, systems determined by locally contracting state maps $F$ always determine a generalized synchronization that has the same degree of regularity as the dynamical system $\phi $, the system $F$, and the  observations $\omega $ that drive it.

\begin{theorem}
\label{main theorem of letter}
Let $\phi \in {\rm Hom} (M) $ be an invertible discrete-time dynamical system on a compact topological space $M $, $\omega \in C ^0(M, \mathbb{R}^d) $ a continuous  observation map, and let $F: \mathbb{R} ^N\times \mathbb{R} ^d \longrightarrow   \mathbb{R} ^N $ be a continuous state map that is locally contracting in the state variable. Let $V \subset \mathbb{R}^N$ be a closed set such that $F(V \times \omega(M))\subset V $ and where $F$ is state-contracting. Then: 
\begin{description}[wide, labelwidth=!, labelindent=0pt]
\item [(i)] The  corresponding restricted system $F:V \times \omega(M) \longrightarrow  V $ has the $(\phi, \omega) $-ESP and hence the generalized synchronization $f_{(\phi, \omega,F)} $ is well-defined. In this case, the GS $f_{(\phi, \omega,F)}: M \longrightarrow V$ is continuous and it is the only map with that codomain that satisfies the identity \eqref{SSM recursion for lemma}.
\item [(ii)] The same conclusion holds when $M$ is a compact differentiable manifold, $\phi \in {\rm Diff}^1(M)$ is an invertible differentiable dynamical system, $\omega $ is of class $C ^1$, $F$ is of class $C ^2 $, and $V$ is closed and convex. In this case, if
\begin{equation}
\label{condition for lfx in differentiable statement}
L_{F _x}< \min \left\{1, 1/ \left\|T \phi ^{-1}\right\|_{\infty}\right\},
\end{equation}
then the associated GS $f_{(\phi, \omega,F)} $ is continuously differentiable.
In \eqref{condition for lfx in differentiable statement} $L_{F _x}:=\sup_{(\mathbf{x}, {\bf z})\in V \times \omega(M)} \left\{ \left\|D _x F(\mathbf{x}, {\bf z})\right\|\right\}$ and $\left\|T \phi\right\|_{\infty}:=\sup_{m \in M} \left\{ \left\|T _m \phi\right\|\right\} $, with $T _m\phi: T _mM \longrightarrow T_{\phi(m)}M $ the tangent map of $\phi$ at $m \in M $.
\end{description}
\end{theorem}

The strategy for obtaining the two parts of this statement consists in using the natural Banach space structures of the spaces $C ^0(M, \mathbb{R} ^N) $ and $C ^1(M, \mathbb{R} ^N) $ of continuous  and  differentiable  functions respectively, to apply the Banach Contraction-Mapping Principle to an automorphism $\Psi: C ^0(M, \mathbb{R} ^N) \longrightarrow C ^0(M, \mathbb{R} ^N) $ of $C ^0(M, \mathbb{R} ^N) $ (respectively of $C ^1(M, \mathbb{R} ^N) $) defined using the right hand side of \eqref{SSM recursion for lemma}, namely, 
\begin{equation*}
\Psi(f) (m):= F \left(f \left(\phi^{-1} (m)\right), \omega (m)\right), \  \mbox{for all} \quad  m \in M.
\end{equation*}
The compactness of $M$ is a crucial hypothesis in all these constructions. Using this approach, the local contractivity hypotheses on $F$ (as well as condition \eqref{condition for lfx in differentiable statement} in the differentiable case) imply that $\Psi $ is  contracting. The generalized synchronization $f_{(\phi, \omega,F)}  $  hence arises as its unique fixed point, a condition that amounts to the identity \eqref{SSM recursion for lemma}.

It is important to emphasize  that for a given state map $F: \mathbb{R}^N \times \omega(M) \longrightarrow \mathbb{R}^N$ there could exist various disjoint closed sets like $V$ that satisfy the hypotheses in the theorem. The use of the restrictions of $F$ to each of them leads in general to {\it different} generalized synchronizations $f_{(\phi, \omega,F)}$ with disjoint codomains. This feature is much related with multistability phenomena and the so-called echo index, as presented in \cite{livi:multistability}. See also \cite{ceni:ashwin:paper1, manjunath:prsl} for related discussions.

The strategy followed in the theorem guarantees that $f_{(\phi, \omega,F)} $ can be obtained by iterating the map $\Psi $ using any continuous or differentiable function  $f _0$ (for instance a constant function) as initial condition. In other words, $f_{(\phi, \omega,F)} $ is the uniform limit of the sequence determined by the iterations:
\begin{equation}
\label{fixed point as limit paper}
f _{n+1}= \Psi(f _n), \quad \mbox{with} \quad f _0:= {\rm constant}.
\end{equation}
There is a fundamental practical difference between the construction of the generalized synchronization using the recursion \eqref{fixed point as limit paper} or via the definition $f_{(\phi, \omega,F)} :=p _0 \circ U ^F\circ S _{(\phi, \omega)} $ in \eqref{definition our gs}. The former requires full knowledge about the dynamical system $\phi$ while the latter only uses its  $\omega $-observations. This difference is of much relevance when a synchronization has to be constructed or {\it learned} using only temporal traces of observations of a given data generating dynamical system.

The existence of differentiable synchronizations for invertible chaotic systems on compact manifolds is not a new fact as it  follows from Takens' Theorem \cite{takensembedding, huke:2006}. Indeed, this result shows that in the presence of certain non-resonance conditions and for generic scalar observations $\omega \in C ^2(M, \mathbb{R} ) $ of a dynamical system $\phi\in {\rm Diff}^1(M) $, with $M$ compact and $q$-dimensional, a $(2q+1)$-truncated version $S _{(\phi, \omega)} ^{2q+1}$ of the the $(\phi, \omega)$-{delay map} given by  $S _{(\phi, \omega)}^{2q+1}(m):=\left(\omega (m), \omega(\phi ^{-1} (m)), \ldots, \omega(\phi ^{-2q} (m))\right)$ is a continuously differentiable embedding. This map is in turn the GS corresponding to the linear state map $F(\mathbf{x}, z):=A \mathbf{x}+ \mathbf{C} z $, with $A$ the lower shift matrix in dimension $2q+1  $ and $\mathbf{C}= (1,0, \ldots,0) \in \mathbb{R}^{2q+1} $ which, by Takens' Theorem, constitutes a differentiable GS for the scalar observations of $\phi $. Theorem \ref{main theorem of letter} allows us to {extend} a part of this statement. Indeed, consider an arbitrary linear state map $F(\mathbf{x}, {\bf z}):=A \mathbf{x}+ \mathbf{C} {\bf z}$, $A \in \mathbb{M} _{N,N} $, $\mathbf{C} \in \mathbb{M} _{N,d} $,  whose connectivity matrix $A$ has singular values bounded by one, that is $\sigma_{{\rm max}}(A)<1 $, and that is driven by continuous $d$-dimensional observations $\omega \in C ^0(M, \mathbb{R}^d) $ of $\phi $.  Theorem \ref{main theorem of letter} guarantees that the system associated to $F$ yields a continuous synchronization $f_{(\phi, \omega,F)} $. Additionally, if $\omega  $ is of class $C ^1 $ (not $C ^2$ as in Takens' Theorem) and $\sigma_{{\rm max}}(A)<\min \left\{1, 1/\left\|T \phi ^{-1}\right\|_{\infty}\right\} $ then $f_{(\phi, \omega,F)} $ is necessarily differentiable. {We note in passing that the encoding of the Takens' delay embedding as a GS shows that the local contractivity invoked in the hypotheses of Theorem \ref{main theorem of letter} is sufficient but not necessary. Moreover, this result guarantees continuity/differentiability for GSs but does not warrant injectivity or immersivity for which, as we know from the Takens' Theorem case, it is likely to be required additional genericity and/or dimensional conditions.}

More recently a similar fact has been proved in \cite[Theorem 2.2.2]{hart:ESNs} for systems defined by Echo State Networks, that is, recurrent neural network-like state systems of the form $F(\mathbf{x}, {\bf z}):=\boldsymbol{\sigma} \left(A \mathbf{x}+ \mathbf{C} {\bf z} + \boldsymbol{\zeta} \right)$, where $A \in \mathbb{M} _{N,N} $, $\mathbf{C} \in \mathbb{M} _{N,d}$, $\boldsymbol{\zeta} \in \mathbb{R}^N $, and the function $\boldsymbol{\sigma}: \mathbb{R}^N \longrightarrow \mathbb{R}^N  $ is constructed by componentwise application of a continuous squashing function $\sigma:\mathbb{R} \longrightarrow \mathbb{R}$. Theorem \ref{main theorem of letter} guarantees that when $F$ is contracting in the state variable then  it  yields a continuous synchronization $f_{(\phi, \omega,F)} $. Additionally, if $\omega   \in C ^1(M, \mathbb{R} ^d)$, $\sigma \in C ^2(\mathbb{R}) $, and $\sigma_{{\rm max}}(A)L _\sigma<\min \left\{1, 1/\left\|T \phi ^{-1}\right\|_{\infty}\right\} $ with $L _\sigma:=\sup_{z \in \mathbb{R}}\left\{| \sigma'(z)|\right\}$, then $f_{(\phi, \omega,F)} $ is necessarily continuously differentiable. This condition coincides with the one in \cite[Theorem 2.2.2]{hart:ESNs}.

We emphasize that the differentiability results that we have established in this paper obviously do not imply that the generalized synchronizations that we introduced are diffeormorphisms onto their images {or even injective. As it was explained in the Introduction and in Section \ref{Definitions and preliminary discussion} such feature is very important at the time of using these results to, for instance, learn attractors from time series.}  It is clear that additional conditions (among others dimensional) need to be required for this to hold. Even though conjectures in this direction have already been formulated (see the last paragraphs in \cite{kocarev1995general}) this question remains open to our knowledge and it will be the subject of a forthcoming publication that is now in preparation \cite{RC19}.

\section{\label{Numerical Illustration}Numerical Illustration}

\begin{figure}
  \caption{A trajectory of the Lorenz system.}
  \label{Lorenz_system}
  \centering
    \includegraphics[width=0.45\textwidth]{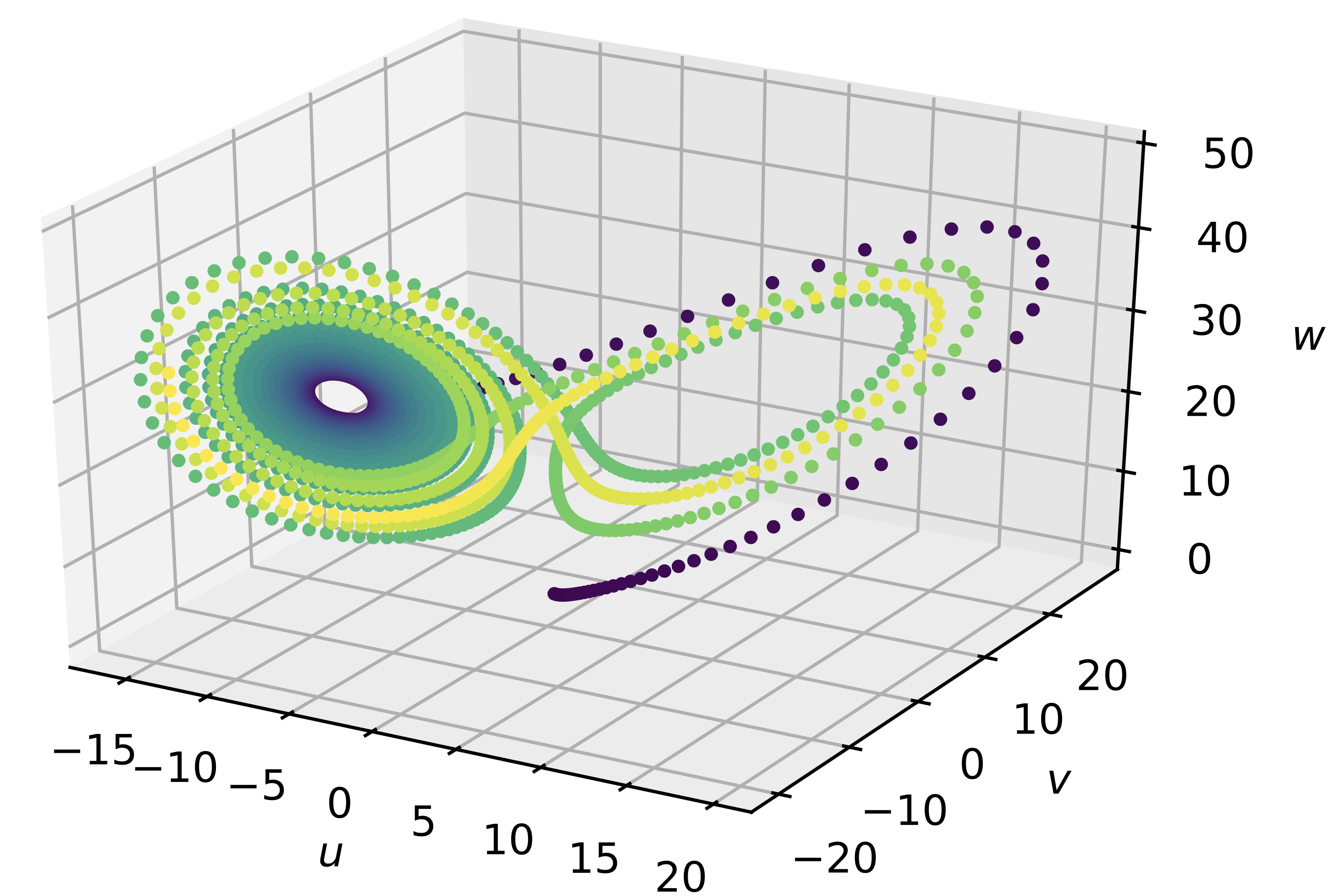}
\end{figure}

\begin{figure}
  \caption{$u$-component of a trajectory of the Lorenz system.}
  \label{x_trajectory}
  \centering
    \includegraphics[width=0.45\textwidth]{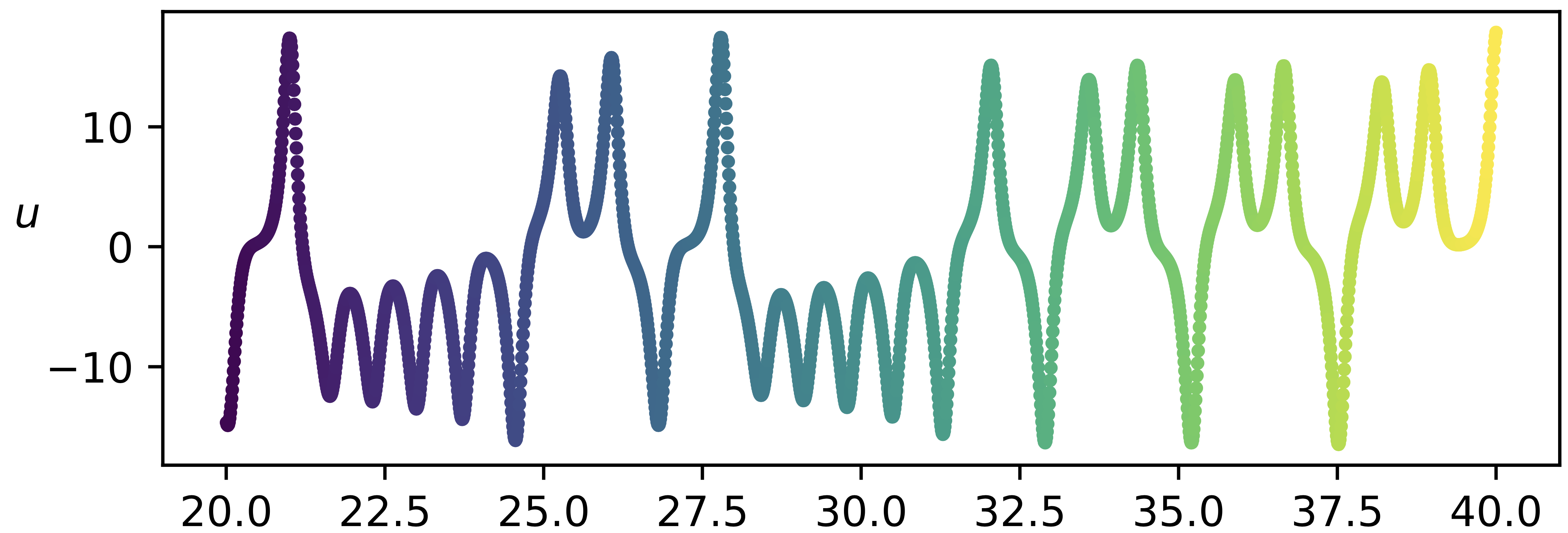}
\end{figure}

In this Section we illustrate various GSs for the Lorenz system using a state map that satisfies the hypotheses of our main result. We recall that the Lorenz system with the parameter values given in the original paper \cite{lorenz1963deterministic} is determined by the differential equation
\begin{align*}
    \dot{u} &= 10(u-v) \\
    \dot{v} &= u(28-w) - v \\
    \dot{w} &= uv - 8w/3.
\end{align*}
A discrete-time dynamical system $\phi$ can be derived from the Lorenz equation like so:
\begin{align*}
    \phi(u,v,w) = (u,v,w) + \int_0^h (\dot{u},\dot{v},\dot{w}) \ dt,
\end{align*}
with $u(0) = u$, $v(0) = v$, $w(0) = w$ and time-step $h$. We simulated a 4000 point (40 time units) trajectory originating from the initial point $m = (0,1,1.05)$ with time-step  $h = 0.01$. Figure \ref{Lorenz_system} shows this trajectory for times $t \in (20,40) $. 
{If we observe only the $u$-component of this trajectory, then the observation function is $\omega(u,v,w) = u$. The corresponding observed trajectory is illustrated in Figure \ref{x_trajectory}.}

\begin{figure}
  \caption{
A phase portrait of the autonomous system $F(\mathbf{x}, z) = (x _1^{\alpha}, x _2^{\alpha}, x _3^{\alpha})$ at the cross section $x _3 = 1$ with the 4 stable fixed points.}
  \label{phase_portrait}
    \includegraphics[width=0.5\textwidth]{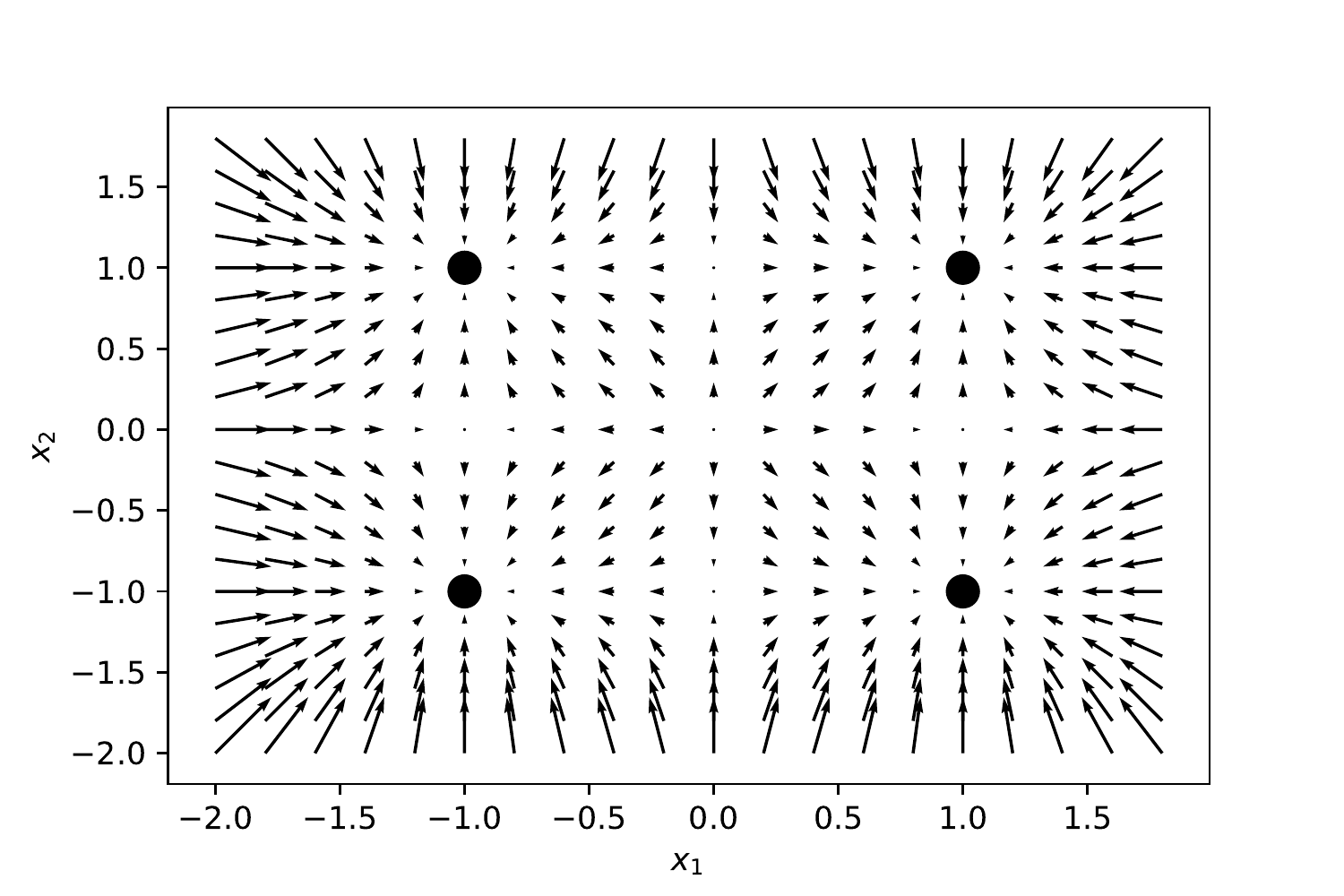}
\end{figure}

We now define the state system  $F : \mathbb{R}^3 \times \mathbb{R} \to \mathbb{R}^3$ given by
\begin{align*}
    F(\mathbf{x}, z) = (x _1^{\alpha},x _2^{\alpha},x _3^{\alpha}) + \lambda(\sin(k z),\cos(k z), \sin^2(k z))
\end{align*}
with $\lambda, k \geq 0$ and $\alpha \in (0,1)$. If we choose $\lambda = 0$ then the state system is autonomous and has 8 stable fixed points at $(1,1,1)$,
$(-1,1,1)$, $(1,-1,1)$, $(1,1,-1)$, $(-1,-1,1)$, $(-1,1,-1)$, $(1,-1,-1)$, $(-1,-1,-1)$. A cross section of the phase portrait of this autonomous system at the $w = 1$ plane is shown in Figure \ref{phase_portrait}.

If we now choose $\lambda = 0.009$, $k = 0.1$ and $\alpha = 0.9$, then we can construct a box $V_1 = [0.9,1.1] \times [0.9,1.1] \times [0.9,1.1]$ containing the fixed point $(1,1,1)$. We can create a similar box about each of the other seven fixed points and denote them $V_2, V_3, \ldots, V_8$. We have by construction for each $i = 1, \ldots , 8$ that
$
    F(V_i,\omega(\mathbb{R}^3)) \subset V_i
$
and that $F$ is state contracting on each box $V_i$. Thus, by Theorem \ref{main theorem of letter} the corresponding restricted state maps $F:V_i \times \omega(M) \longrightarrow  V_i $ have the $(\phi, \omega) $-ESP, so the generalized synchronizations $f_{(\phi, \omega,F)} $ are well-defined and differentiable.

\begin{figure}
\centering
  \caption{Image of the Lorenz solution by two different synchronization maps $f_{(\phi, \omega,F)} $ that contain the points $ (\pm 1,1,1)$  in their images.}
  \label{reservoir_lorenz}
\includegraphics[width=0.45\textwidth]{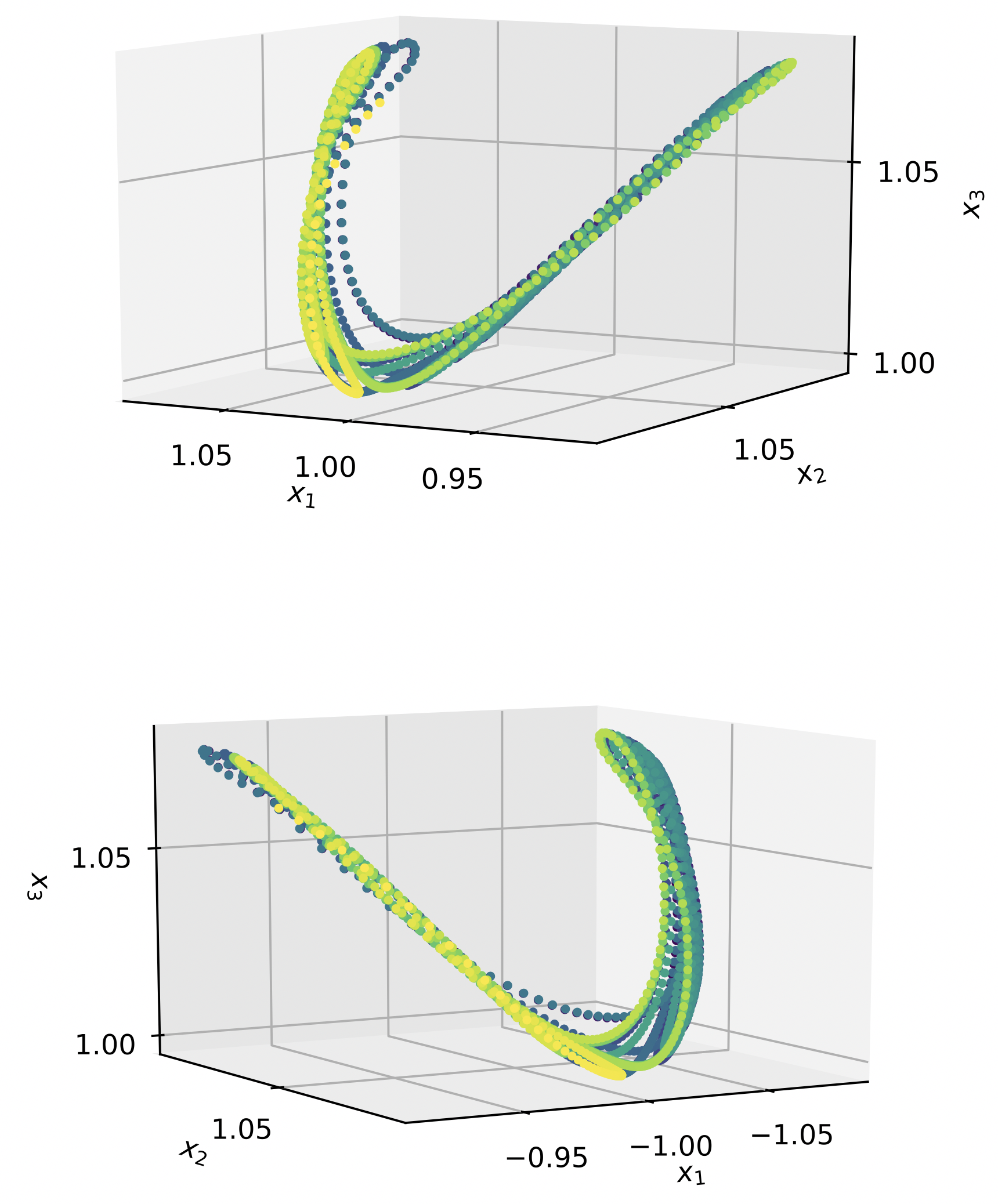}
\end{figure}

We then computed the states
\begin{align*}
    \mathbf{x} = F(\mathbf{x}_{t-1},\omega(\phi^t(m)))
\end{align*}
from two different initial states $\mathbf{x}_0 = (1,1,1)$ and $\mathbf{x}_0 = (-1,1,1)$ and plotted the results in Figure \ref{reservoir_lorenz}. Since $F$ has the fading memory property and hence the input-forgetting property  (see Proposition \ref{contracting implies omega ESP}), the values $\mathbf{x} _t  $  produced by these iterations converge (after a washout period) to the unique state-space solution $\left\{\mathbf{x} _t\right\}_{t \in \Bbb Z}$ determined by the infinite input $\left\{\omega(\phi^t(m))\right\}_{t \in \Bbb Z}$.  We chose the interval $(0,20)$ as washout period. The defining property \eqref{generalized synchronization condition} of the synchronization maps $f_{(\phi, \omega,F)} $ guarantee that, after the washout period, the values $\mathbf{x} _t  $ are virtually indistinguishable from the image $f_{(\phi, \omega,F)} (\phi ^t(m))$. Furthermore, for \emph{any} initial point $x_0 \in V_i$, the states $\mathbf{x}_t$ converge to the \emph{same} image $f_{(\phi, \omega,F)} (\phi ^t(m))$ contained by the box $V_i$.

The mapping $f_{(\phi,\omega,F)}$ from the Lorenz attractor to the reservoir space can be numerically computed using two different methods that yield similar results. We have just described the first, and the second uses the iterations \eqref{fixed point as limit paper}, taking as an initial conditions the constant map $f _0= \mathbf{x} _0 $, with $\mathbf{x} _0 =(1,1,1)$ or $\mathbf{x} _0 =(-1,1,1)$, respectively.

Figure \ref{reservoir_lorenz} shows that that the GS $f_{(\phi, \omega,F)}$ is visually a differentiable mapping of the Lorenz trajectory into the state space. We stress nevertheless that we have no theoretical guarantees (yet) of the injectivity of $f_{(\phi, \omega,F)} (\phi ^t(m))$ and hence of the supervised learnability of the Lorenz attractor using this particular synchronization.

%

\appendix
\section{Proof of Proposition II.1} 

\noindent The proof of  these claims hinges on various results that are already available in the literature. 

\medskip

\noindent {\bf (i)} We start by using an argument similar to the one in \cite[Remark 2]{RC10}. Fix an element $\mathbf{v} \in D _N $. Given that $ \left\{ {\bf v}\right\} \times D_d $ is a compact  subset of  $\mathbb{R} ^N\times \mathbb{R} ^d $ and $F$ is continuous, then (see \cite[Theorem 26.5]{Munkres:topology}) $F( \left\{ {\bf v}\right\} \times D_d) $ is a compact subset of $\mathbb{R}^N $ and hence bounded (see \cite[Theorem 27.3]{Munkres:topology}). This implies the existence of a constant $r>0 $ such that $F( \left\{ {\bf v}\right\} \times D_d) \subset \overline{B _{r}(\mathbf{v})}\cap D _N$. We now see that for any constant $M>r/(1-c)>0 $, the compact set $W:=D _N\cap  \overline{B _{M}(\mathbf{v})} $ satisfies that $F(W\times D_d) \subset W$. Indeed, for any $(\mathbf{x}, \mathbf{z}) \in W \times D_d $ we can write: 
\begin{multline*}
\left\|F(\mathbf{x}, {\bf z})\right\|\leq \left\|F(\mathbf{x}, {\bf z})- F({\bf v}, {\bf z})\right\|+ \left\|F({\bf v}, {\bf z})\right\|\\
\leq c \left\|\mathbf{x}-\mathbf{v}\right\|+ r\leq cM+r< M,
\end{multline*}
which shows that $F(\mathbf{x}, {\bf z}) \in  \overline{B _{M}(\mathbf{v})} $. Since by construction $W \subset D _N $ and by hypothesis $F(D _N\times D_d) \subset D _N$ we also have that $F(\mathbf{x}, {\bf z}) \in D _N $ and hence that $F(\mathbf{x}, {\bf z}) \in W$. Notice that if $D_N $ is convex then so is $W:=D _N\cap  \overline{B _{M}(\mathbf{v})} $.

The $(D _d^{\mathbb{Z}_{-}}, W ^{\mathbb{Z}_{-}})$-ESP of  $F: W\times D_d \longrightarrow   W$ and part {\bf (ii)} can be obtained as in the second part in Theorem 3.1 in \cite{RC7}. The statement on the continuity with respect to the product topology and the FMP with respect to any weighting sequence is a consequence of Corollary 2.7 in the same reference. Similar results under more or less general hypotheses can also be found in Theorem 7 and Theorem 12 of \cite{RC9}. \quad $\blacksquare$

\section{Proof of Lemma II.2} 

\noindent We first recall that, by definition, as $F$ has the $(\phi, \omega) $-ESP, there exists a unique filter $U ^F: S _{(\phi, \omega)}(M) \longrightarrow (D_N)^{\Bbb Z} $ that associates to each orbit $S _{(\phi, \omega)}(m) $ the unique solution sequence $U^F(S _{(\phi, \omega)}(m)) \in (D_N)^{\Bbb Z} $ that satisfies
\begin{equation}
\label{recursion for UF}
U^F(S _{(\phi, \omega)}(m)) _t= F \left(U^F(S _{(\phi, \omega)}(m)) _{t-1}, S _{(\phi, \omega)}(m) _t\right), 
\end{equation}
for all $t \in \Bbb Z $.
The $(\phi, \omega)$-ESP automatically ensures that $U ^F $ is causal and time-invariant (see \cite[Proposition 2.1]{RC7} for a detailed proof) and hence
\begin{equation}
\label{time invariance explicit}
U ^F \circ T _t= T _t \circ U^F, \quad \mbox{for all} \quad t \in \Bbb Z,
\end{equation}
with $T _t $ the time delay operator introduced above \eqref{relation p and T}. With these elements it is now easy to prove \eqref{defining claim for GS}. Indeed, 
\begin{multline}
U^F(S _{(\phi, \omega)}(m)) _t= p _t \left(U^F(S _{(\phi, \omega)}(m)) \right)\\
= p _0 \circ T _{-t} \left(U^F(S _{(\phi, \omega)}(m)) \right)
= p _0  \left(U^F\left( T_{-t} \left(S _{(\phi, \omega)}(m)\right) \right)\right)\\
=p _0 \left(
U^F \left(S _{(\phi, \omega)}(\phi^t(m))\right)\right)=f_{(\phi, \omega,F)} \left(\phi ^t (m)\right),
\end{multline}
as required. In this chain of equalities, we have used \eqref{relation p and T} in the second equality, the time invariance \eqref{time invariance explicit} of $U^F $ in the third one, and the time equivariance \eqref{somega equivariant} of $S _{(\phi, \omega)} $ in the fourth one. \quad $\blacksquare$

\section{Proof of Lemma II.3.} 

\noindent This is a straightforward corollary of the relation \eqref{defining claim for GS} and the recursion \eqref{recursion for UF} satisfied by $U^F $. Indeed, if we set $t = 0  $ in \eqref{defining claim for GS}, we have 
\begin{multline*}
f_{(\phi, \omega,F)} (m)=U ^F(S _{(\phi, \omega)} (m)) _0\\
= F \left(U ^F(S _{(\phi, \omega)} (m)) _{-1}, S _{(\phi, \omega)} (m) _0\right)\\
=F \left(f_{(\phi, \omega,F)} \left(\phi ^{-1} (m)\right), \omega(m)\right),
\end{multline*}
where in the second equality we used \eqref{recursion for UF} and in the third one \eqref{defining claim for GS} again with $t=-1 $.
 \quad $\blacksquare$
 
 \section{Existence and uniqueness of continuous generalized synchronizations}
 
 \begin{proposition}
\label{continuous synchronization}
Let $\phi \in {\rm Hom} (M) $ be an invertible dynamical system on a compact topological space $M $, $\omega \in C ^0(M, \mathbb{R}^d) $ an observation map, and let $F: D _N\times \omega(M)\longrightarrow   D_N $ a continuous state map that is a contraction on the state variable, with $D _N\subset \mathbb{R} ^N $ a closed set. Then: 
\begin{description}
\item [(i)] There exists a compact subset $W \subset D _N $ such that  $F(W \times \omega(M))\subset W $ and the system determined by $F:W \times \omega(M) \longrightarrow W$ and driven by the $\omega $-observations of $\phi $ has the $(\phi, \omega) $-ESP and hence the generalized synchronization $f_{(\phi, \omega,F)}:M \longrightarrow W$ is well-defined. If $D _N $ is convex then $W $ can be chosen to be convex.
\item [(ii)] The map $f_{(\phi, \omega,F)} :M \longrightarrow W$ is continuous and it is the only one with that codomain that satisfies the identity \eqref{SSM recursion for lemma}.
\end{description}
\end{proposition}

\noindent\textbf{Proof.\ \ (i)} This claim is a  consequence of Proposition \ref{contracting implies omega ESP}. Indeed, since $\omega$ is continuous and $M$ is compact then so is $\omega (M)$. The contractivity hypothesis and the closeness of $D_N $ imply by  Proposition \ref{contracting implies omega ESP} the existence of the compact set $W$ in the statement and that the corresponding system $F:W \times \omega(M) \longrightarrow W$ has the $(\omega(M)^{\mathbb{Z}_{-}}, W ^{\mathbb{Z}_{-}}) $-ESP. Given that $S _{(\phi, \omega)}^-(M) \subset  \omega(M)^{\mathbb{Z}_{-}} $ this implies that the $(S _{(\phi, \omega)}^-(M), W ^{\mathbb{Z}_{-}})$-ESP also holds. This condition amounts to the $(\phi, \omega) $-ESP and to the existence of a FMP filter $U ^F:S _{(\phi, \omega)}^-(M)\longrightarrow W ^{\mathbb{Z}_{-}}$ which, in passing, shows that the generalized synchronization $f_{(\phi, \omega,F)} :=p _0 \circ  U ^F \circ S _{(\phi, \omega)} :M \longrightarrow W$ is well-defined.

\medskip

\noindent {\bf (ii) } We first show the continuity of $f_{(\phi, \omega,F)}:= p _0 \circ U ^F \circ  S_{(\phi, \omega)}:M \longrightarrow W$ by noticing that it is a composition of  continuous functions. Indeed, using the notation introduced in the previous point, if we endow the sets $\omega(M)^{\mathbb{Z}_{-}}$ and $W ^{\mathbb{Z}_{-}}$ with the product topology, the map $p _0:W ^{\mathbb{Z}_{-}} \longrightarrow W$ is clearly continuous, $U^F: \omega(M)^{\mathbb{Z}_{-}}\longrightarrow W ^{\mathbb{Z}_{-}} $ is continuous by Proposition \ref{contracting implies omega ESP}, and $S_{(\phi, \omega)}: M \longrightarrow \omega(M)^{\mathbb{Z}_{-}}$ is also continuous because it is a Cartesian product of continuous maps (see \cite[Theorem 19.6]{Munkres:topology}).

In order to prove the uniqueness statement, we first endow the set of continuous functions $C ^0(M, \mathbb{R} ^N) $ with the Banach space structure induced by the norm
\begin{equation}
\label{infty norm for continuous functions}
\left\|f\right\|_{\infty}:=\sup_{m \in M} \left\{\left\|f(m)\right\|\right\}, \quad f \in C ^0(M, \mathbb{R} ^N).
\end{equation}
We emphasize that this norm is always finite due to the compactness of $M$. Moreover, it is easy to show that the subset $C ^0(M, W)$ of $C ^0(M, \mathbb{R} ^N) $ made of functions that have the compact subset $W \subset \mathbb{R} ^N $ as codomain is closed with respect to the topology generated by the norm $\left\|\cdot \right\|_{\infty}$. Indeed, let $f \in C ^0(M, \mathbb{R} ^N) $ be an element in the closure of $C ^0(M, W)$ and let $\left\{f _n\right\}_{n \in \mathbb{N}} $ be a sequence of elements in $C ^0(M, W)$ that have $f$ as limit. This implies that for any $\epsilon>0 $ there exists $N(\epsilon) \in \mathbb{N} $ such that for any $n >N(\epsilon) $ we have that $\left\|f _n-f \right\|_{\infty}< \epsilon $. This implies, in particular, that for any point $ m \in M $:
\begin{equation*}
\left\|f _n (m)- f (m)\right\|\leq \left\|f _n-f \right\|_{\infty}< \epsilon 
\end{equation*}
and hence $\lim\limits_{n \rightarrow \infty} f _n(m)= f (m) $, which shows that $f (m) $ is a limit point of $W$ and hence belongs to $W$ as it is a closed set. Since $m \in M $ is arbitrary this shows that $f \in C ^0(M, W)$ and hence $C ^0(M, W)$  is closed in $C ^0(M, \mathbb{R} ^N) $.

We now shall prove the uniqueness statement by showing that the map $\Psi:C ^0(M, W) \longrightarrow C ^0(M, W) $ defined by
\begin{equation}
\label{definition of Psi b}
\Psi(f) (m):= F \left(f \left(\phi^{-1} (m)\right), \omega (m)\right),
\end{equation}
for all $m \in M,\, f \in C ^0(M, \mathbb{R} ^N) $,
is well-defined and a contraction with respect to the norm introduced in \eqref{infty norm for continuous functions}. The result follows then from the Banach Contraction-Mapping Principle (see \cite[Theorem 3.2]{Shapiro:Farrago}).

First, the continuity hypotheses on $\phi, \omega$, and $F:W \times \omega(M) \longrightarrow W$ imply that $\Psi(f) \in C ^0(M, W) $ whenever $f \in C ^0(M, W) $. Now, let $f,g \in C ^0(M, W) $. Then, since we assume that $F$  is a contraction in the state variable with constant $0<c<1$, we have that,
\begin{multline}
\label{phi contraction proof}
\left\|\Psi(f)-\Psi(g)\right\|_{\infty}\\
=\sup _{m \in M} \{\|F \left(f \left(\phi^{-1} (m)\right), \omega (m)\right)\\
-F \left(g \left(\phi^{-1} (m)\right), \omega (m)\right)   \|\}\\
\leq c \sup _{m \in M} \left\{\left\|f \left(\phi^{-1} (m)\right)-g \left(\phi^{-1} (m)\right)   \right\|\right\}=c \left\|f-g\right\|_{\infty},
\end{multline}
which shows that $\Psi  $ is a contraction with respect to the norm in \eqref{infty norm for continuous functions} and hence there exists a unique element in $C ^0(M, W) $ that satisfies the identity \eqref{SSM recursion for lemma}. As we just proved that $f_{(\phi, \omega,F)} \in  C ^0(M, W)$ and by Lemma \ref{lemma SSM recursion} the state synchronization map $f_{(\phi, \omega,F)} $  satisfies  \eqref{SSM recursion for lemma}, that unique element is $f_{(\phi, \omega,F)} $ necessarily. \quad $\blacksquare$ 

\begin{remark}
\label{global vs lokal 1}
\normalfont
It is worth emphasizing  that, for a given state map $F: \mathbb{R}^N \times \omega(M) \longrightarrow \mathbb{R}^N$, there could exist various disjoint closed sets like $D_N$ that satisfy the hypotheses of the proposition. The use of the restrictions $F: D _N \times \omega(M) \longrightarrow D_N$ to each of them leads in general to {\it different} generalized synchronizations $f_{(\phi, \omega,F)}$ whose codomains are compact subsets of each of the closed set $D_N $. This feature is much related with multistability phenomena and the so-called echo index, as presented in \cite{livi:multistability}. See also \cite{ceni:ashwin:paper1, manjunath:prsl} for related discussions.
\end{remark}

\begin{remark}
\normalfont
Having obtained in the proof the state synchronization map $f_{(\phi, \omega,F)}:M \longrightarrow W $ as the unique fixed point of a contracting map, a standard result about maps of this type on metric spaces (see \cite[Proposition 3.4]{Shapiro:Farrago}), guarantees that $f_{(\phi, \omega,F)} $ can be obtained by iterating the map $\Psi $ defined in \eqref{definition of Psi b} using any function $f _0 \in C ^0(M, W)$ (for instance a constant function) as initial condition. In other words, $f_{(\phi, \omega,F)} $ is the uniform limit of the sequence of functions determined by the iterations:
\begin{equation}
\label{fixed point as limit}
f _{n+1}= \Psi(f _n), \quad \mbox{with} \quad f _0:= \mathbf{w} \in W,
\end{equation}
where the constant element $\mathbf{w} \in W $ is arbitrary.
\end{remark}

\begin{remark}
\normalfont
There is a fundamental difference of much practical importance between the construction of the generalized synchronization $f_{(\phi, \omega,F)} $ via the iteration of the map $\Psi $ as in \eqref{fixed point as limit} and using the definition $f_{(\phi, \omega,F)} :=p _0 \circ U ^F\circ S _{(\phi, \omega)} $. The former requires full knowledge about the dynamical system $\phi$ while the latter only uses its  $\omega $-observations. This difference is of much relevance when a synchronization has to be constructed or {\it learned} using only observations of a given data generating dynamical system.
\end{remark}

\section{Existence and uniqueness of differentiable synchronizations}

In this section we extend the result in Proposition \ref{continuous synchronization} and we show that when the dynamical system $\phi \in {\rm Diff}^1(M)$ and the maps $\omega $  and $F$  are differentiable then the state synchronization map $f_{(\phi, \omega,F)} $ is necessarily differentiable. The strategy that we follow also consists in proving that the map $\Psi $ introduced in \eqref{definition of Psi b} is a contraction but, this time around, on a well-chosen closed subset of $C ^1(M, \mathbb{R} ^N)$. We start with some preliminaries.

\paragraph{Banach space structures in $C ^1(M, \mathbb{R} ^N)$.} All along this section we assume that $M$ is a compact, connected, Hausdorff, and second-countable manifold and hence it can be endowed with a Riemannian metric $g$  (see \cite[Proposition 2.10]{do:carmo:1993}). Now, for any $f \in C ^1(M, \mathbb{R} ^N)$, define
\begin{equation*}
\left\|Df\right\|_{\infty}=\sup_{m \in M} \left\{ \left\|Df(m)\right\|\right\} \quad \mbox{with}
\end{equation*}
\begin{equation*}
\left\|Df(m)\right\|=\sup_{\stackanchor{\scriptstyle \mathbf{v} \in T _mM}{\scriptstyle \mathbf{v}\neq {\bf 0}} } 
\left\{\frac{\left\|Df(m) \cdot \mathbf{v}\right\|}{\left(g(m)(\mathbf{v},\mathbf{v})\right)^{1/2}}\right\}.
\end{equation*}
Analogously, if $\phi:M \rightarrow M $ is a $C ^1 $ map, we can define:
\begin{equation*}
\left\|T \phi\right\|_{\infty}=\sup_{m \in M} \left\{ \left\|T _m \phi\right\|\right\} \quad \mbox{with} 
\end{equation*}
\begin{equation*}
\left\|T _m \phi\right\|=\sup_{\stackanchor{\scriptstyle \mathbf{v} \in T _mM}{\scriptstyle \mathbf{v}\neq {\bf 0}} } 
\left\{\frac{\left(g(\phi(m))(T _m\phi \cdot \mathbf{v},T _m\phi \cdot \mathbf{v})\right)^{1/2}}{\left(g(m)(\mathbf{v},\mathbf{v})\right)^{1/2}}\right\}.
\end{equation*}

It can be proved by using the results in Chapter 2 of \cite{abraham:robbin} that, for any $\delta>0 $, the norms $\left\|\cdot \right\|_{C ^1(\delta)}$ defined by
\begin{equation}
\label{definition c1d}
\left\|f \right\|_{C ^1(\delta)}:= \left\|f\right\|_{\infty}+ \delta\left\|Df\right\|_{\infty}
\end{equation}
endow $C ^1(M, \mathbb{R} ^N)$ with a Banach space structure. It can also be shown (see \cite[Theorem 11.2 {\bf (ii)}]{abraham:robbin}) that all these norms generate the same topology in $C ^1(M, \mathbb{R} ^N)$ that is also independent of the choice of Riemannian metric $g$ and coincides with the weak and strong topologies introduced in Chapter~2 of \cite{Hirsch:book}. 

We introduce two technical lemmas that are needed later on.

\begin{lemma}
\label{c1d closed}
In the setup that we just described, let $R>0 $, $W \subset \mathbb{R}^N  $ a closed set, and define:
\begin{equation*}
\Omega(R,W):= \left\{f \in C ^1(M, W)\mid \left\|Df\right\| _{\infty}\leq R\right\},
\end{equation*}
where $f \in C ^1(M, W) $ is by definition the subset of $f \in C ^1(M, \mathbb{R} ^N) $ made by the functions which take values in $W \subset \mathbb{R} ^N $.
The set $\Omega (R,W) $ is closed in $\left( C ^1(M, \mathbb{R} ^N), \left\|\cdot \right\|_{C ^1(\delta)} \right)$ for any $\delta>0 $.
\end{lemma}

\noindent\textbf{Proof.\ \ }We proceed by showing that $\Omega(R,W) $ is the intersection of two closed sets, that is, $\Omega(R,W)=\Omega(R)\cap  C ^1(M, W) $, where $\Omega(R):=  \left\{f \in C ^1(M, \mathbb{R}^N)\mid \left\|Df\right\| _{\infty}\leq R\right\}$. 

We first prove that $\Omega(R)$ is closed by showing that the complementary set $\Omega (R) ^c$ is open in $\left( C ^1(M, \mathbb{R} ^N), \left\|\cdot \right\|_{C ^1(\delta)} \right)$. Let $f \in \Omega (R) ^c $ such that $\left\|Df\right\|_{\infty}=K>R $. Given $\epsilon:= \delta(K-R)$ we will show that the ball $B_{C ^1(\delta)}(f, \epsilon) $  is a subset of $\Omega (R) ^c$. Indeed, for any $g \in  B_{C ^1(\delta)}(f, \epsilon) $:
\begin{multline*}
K= \left\|Df\right\|_{\infty}=\left\|Df+Dg-Dg\right\|_{\infty}\\
\leq \left\|Df-Dg\right\| _{\infty}+ \left\|Dg\right\| _{\infty}\\
\leq \frac{1}{\delta} \left\|f-g\right\|_{\infty}+\left\|Df-Dg\right\| _{\infty}+ \left\|Dg\right\| _{\infty}\\
= \frac{1}{\delta} \left\|f-g\right\|_{C ^1(\delta)}+ \left\|D g\right\|_{\infty}< \frac{\epsilon}{\delta}+\left\|D g\right\|_{\infty}\\
=K-R+\left\|D g\right\|_{\infty},
\end{multline*}
which implies that $\left\|D g\right\|_{\infty}>R $ and hence that $B_{C ^1(\delta)}(f, \epsilon) \subset \Omega (R) ^c $. 

We now show that $C ^1(M, W) $ is closed in $C ^1(M, \mathbb{R} ^N)$. Let $f \in C ^1(M, \mathbb{R} ^N) $ be an element in the closure of $C ^1(M, W)$ and let $\left\{f _n\right\}_{n \in \mathbb{N}} $ be a sequence of elements in $C ^1(M, W)$ that have $f$ as limit. This implies that for any $\epsilon>0 $ there exists $N(\epsilon) \in \mathbb{N} $ such that for any $n >N(\epsilon) $ we have that $\left\|f _n-f \right\|_{C ^1(\delta)}< \epsilon $. This implies, in particular, that for any point $ m \in M $:
\begin{multline*}
\left\|f _n (m)- f (m)\right\|\leq \left\|f _n-f \right\|_{\infty}+ \delta\left\|Df _n-Df \right\|_{\infty}\\
= \left\|f _n-f\right\|_{C ^1(\delta)}< \epsilon 
\end{multline*}
and hence $\lim\limits_{n \rightarrow \infty} f _n(m)= f (m) $, which shows that $f (m) $ is a limit point of $W$ and hence belongs to $W$ as it is a closed set. Since $m \in M $ is arbitrary this shows that $f \in C ^1(M, W)$ and hence $C ^1(M, W)$  is closed in $C ^1(M, \mathbb{R} ^N) $. \quad $\blacksquare$

\medskip

The following notation will be used in the sequel. Let $D _N\subset \mathbb{R} ^N  $ and $D _d \subset {\Bbb R}^d $ be open subsets and  $F\in C ^2(D_N \times D_d, D_N) $. The symbols $D _x F(\mathbf{x}, {\bf z}) $, $D _z F(\mathbf{x}, {\bf z}) $, $D _{xx} F(\mathbf{x}, {\bf z}) $, and $D _{xz} F(\mathbf{x}, {\bf z}) $ denote the partial derivatives of $F$ with respect to the variables indicated in the subindices at the point $(\mathbf{x}, {\bf z}) \in D_N \times D_d$.

\begin{lemma}
\label{psi restricts to omeagar}
Let $\phi \in {\rm Diff}^1(M)$ be a dynamical system on the compact differentiable manifold $M$ and consider the observation $\omega \in C ^1(M, \mathbb{R}^d) $  and state  $F\in C ^1(D_N \times D_d, D_N) $  maps, with $D _N\subset \mathbb{R} ^N  $ and $D _d \subset {\Bbb R}^d $ open subsets. Let $W\subset D_N$ be a subset and suppose that $\omega(M) \subset D_d$ and that $F(W \times \omega(M) )\subset W$. Define,
\begin{equation}
\label{Ls for F 1}
\begin{array}{lcl}
L_{F _x}&:=&\sup_{(\mathbf{x}, {\bf z})\in W \times \omega(M)} \left\{ \left\|D _x F(\mathbf{x}, {\bf z})\right\|\right\}, \\
 L_{F _z}&:= &\sup_{(\mathbf{x}, {\bf z})\in W \times \omega(M)} \left\{ \left\|D _z F(\mathbf{x}, {\bf z})\right\|\right\}.
 \end{array}
\end{equation}
Suppose that the constants $L_{F _x} $ and $L_{F _z}$ defined in \eqref{Ls for F 1} are finite and that $L_{F _x}\left\|T \phi^{-1}\right\|_{\infty}<1 $. Choose a constant $R>0$ such that 
\begin{equation}
\label{hypothesis on R}
R> \frac{L_{F _z} \left\|D \omega\right\|_{\infty}}{1-L_{F _x}\left\|T \phi^{-1}\right\|_{\infty}}.
\end{equation}
Then, the map $\Psi  $ introduced in \eqref{definition of Psi b} maps the space $C ^1(M, W) $ into itself and, additionally, it restricts to $\Omega(R,W) $, that is, $\Psi(\Omega(R,W)) \subset \Omega (R,W) $.
\end{lemma}

\noindent\textbf{Proof.\ \ } We first note that, for any $ f \in C ^1(M, \mathbb{R} ^N) $, its image $\Psi(f) $ belongs necessarily to  $C ^1(M, \mathbb{R} ^N) $ as it is a composition of $C ^1  $ maps. Moreover, the hypothesis $F(W \times \omega(M) )\subset W)$ guarantees that $\Psi $ preserves $C ^1(M, W) $.

We now prove the inclusion $\Psi(\Omega(R,W)) \subset \Omega (R,W) $. Let $f \in \Omega (R,W) $. By the chain rule and the definitions \eqref{Ls for F 1}, for any $m \in M $  and $\mathbf{v}\in T _m M $:
\begin{multline*}
\left\|D \Psi(f)(m) \cdot \mathbf{v}\right\|\\
= \|   D _xF(f(\phi^{-1}(m)), \omega (m))\cdot D f (\phi^{-1} (m))\cdot T _m \phi^{-1} \cdot \mathbf{v}\\+
D _zF(f(\phi^{-1}(m)), \omega (m))\cdot D\omega (m) \cdot \mathbf{v}\|\\
\leq
\left(L_{F _x} \left\|D f\right\|_{\infty}\left\|T \phi ^{-1}\right\|_{\infty}+ L_{F _z}\left\|D \omega\right\|_{\infty}\right) \left\|\mathbf{v}\right\|\\
\leq \left(L_{F _x} R\left\|T \phi ^{-1}\right\|_{\infty}+ L_{F _z}\left\|D \omega\right\|_{\infty}\right) \left\|\mathbf{v}\right\|,
\end{multline*}
which, using the hypothesis \eqref{hypothesis on R}, implies that
\begin{equation*}
\left\|D\Psi( f) \right\|_{\infty}\leq L_{F _x} R\left\|T \phi ^{-1}\right\|_{\infty}+ L_{F _z}\left\|D \omega\right\|_{\infty} < R,
\end{equation*}
and guarantees that $\Psi(f) \in \Omega (R,W) $, as required. \quad $\blacksquare$ 

\paragraph{The main theorem.} The following result extends Proposition \ref{continuous synchronization} to the differentiable case. The statement requires an extension of the constants introduced in \eqref{Ls for F 1} to higher order derivatives and with respect to a set $V \subset D_N $, that is:
\begin{equation}
\label{Ls for F}
\begin{array}{lcl}
L_{F _x}&:=&\sup_{(\mathbf{x}, {\bf z})\in V \times \omega(M)} \left\{ \left\|D _x F(\mathbf{x}, {\bf z})\right\|\right\},\\
L_{F _z}&:=&\sup_{(\mathbf{x}, {\bf z})\in V \times \omega(M)} \left\{ \left\|D _z F(\mathbf{x}, {\bf z})\right\|\right\},\\
L_{F _{xx}}&:=&\sup_{(\mathbf{x}, {\bf z})\in V \times \omega(M)} \left\{ \left\|D _{xx} F(\mathbf{x}, {\bf z})\right\|\right\},\\
L_{F _{xz}}&:=&\sup_{(\mathbf{x}, {\bf z})\in V \times \omega(M)} \left\{ \left\|D _{xz} F(\mathbf{x}, {\bf z})\right\|\right\}.
\end{array}
\end{equation}

\begin{theorem}
\label{differentiable generalized synchronizations theorem}
Let $\phi \in {\rm Diff}^1(M)$ be a dynamical system on the compact manifold $M$ and consider the observation $\omega \in C ^1(M, \mathbb{R}^d) $  and state  $F\in C ^2(D_N \times D_d, D_N) $  maps, with $D _N\subset \mathbb{R} ^N  $ and $D _d \subset {\Bbb R}^d $ open subsets such that $\omega(M) \subset D_d$. Let $V\subset D_N$ be a closed convex subset and suppose that  $F(V \times \omega(M) )\subset V$. Suppose that the bounds for the partial derivatives of $F$  introduced in \eqref{Ls for F} are all finite and that, additionally,
\begin{equation}
\label{condition for lfx in differentiable}
L_{F _x}< \min \left\{1, 1/ \left\|T \phi ^{-1}\right\|_{\infty}\right\}.
\end{equation}
Then there exists a compact and convex subset $W \subset V$ such that  $F(W \times \omega(M))\subset W $ and: 
\begin{description}
\item [(i)] The system determined by $F:W \times \omega(M) \longrightarrow W$ and driven by the $\omega $-observations of $\phi $ has the $(\phi, \omega) $-ESP and hence the generalized synchronization $f_{(\phi, \omega,F)}: M \longrightarrow W $ is well-defined.
\item [(ii)] The map $f_{(\phi, \omega,F)} $ belongs to $C ^1(M, W) $ and it is the only one that satisfies the identity \eqref{SSM recursion for lemma}.
\end{description}
\end{theorem}

\noindent\textbf{Proof of the theorem.\ \ (i)} First of all, since by hypothesis $0<L_{F _x}<1 $, then $F:V \times \omega(M) \longrightarrow V$ is a contraction in the state variable. Indeed, by the Mean Value Theorem and the convexity hypothesis on $V$, for any $\mathbf{x} _1, \mathbf{x} _2 \in V$, ${\bf z} \in \omega(M)$:
\begin{multline}
\label{f contraction with lfx}
\left\|F(\mathbf{x} _1, {\bf z})-F(\mathbf{x} _2, {\bf z})\right\|\leq 
\sup_{\mathbf{x} \in V}\left\{\left\|D _xF(\mathbf{x}, {\bf z})\right\|\right\} \left\|\mathbf{x} _1- \mathbf{x} _2\right\|\\
\leq L_{F _x} \left\|\mathbf{x} _1- \mathbf{x} _2\right\|.
\end{multline}
The claim then follows from part {\bf (i)} in Proposition \ref{continuous synchronization}.

\medskip

\noindent {\bf (ii)} The strategy of the proof consists in showing that, in the hypotheses in the statement, and using the restricted state map $F:W \times \omega(M) \longrightarrow W$ obtained in part {\bf (i)}, the map $\Psi: \Omega(R,W) \longrightarrow \Omega(R,W) $ introduced in \eqref{definition of Psi b} and with $R> 0$ satisfying the condition \eqref{hypothesis on R} is a contraction for some norm $\left\|\cdot \right\|_{C ^1(\delta_0)}$ of the type defined in \eqref{definition c1d}. The value $\delta _0>0 $ will be specified later on. 

This construction is feasible because the hypothesis \eqref{condition for lfx in differentiable}  implies that $L_{F _x}\left\|T \phi^{-1}\right\|_{\infty}<1 $. Moreover, by Lemma \ref{psi restricts to omeagar} and the finiteness of the constants introduced in \eqref{Ls for F}, a constant $R>0 $ can be chosen so that the map $\Psi: \Omega(R,W) \longrightarrow \Omega(R,W) $ is well-defined. Moreover, since by Lemma \ref{c1d closed} the set $\Omega (R,W) $ is closed in $\left( C ^1(M, \mathbb{R} ^N), \left\|\cdot \right\|_{C ^1(\delta)} \right)$ for any $\delta>0 $, it is hence a complete metric space to which the Banach Contraction-Mapping Principle  \cite[Theorem~3.2]{Shapiro:Farrago} can be applied. Consequently, if we are able to prove that $\Psi: \Omega(R,W) \longrightarrow \Omega(R,W) $ is a contraction with respect to the metric inherited from $\left( C ^1(M, \mathbb{R} ^N), \left\|\cdot \right\|_{C ^1(\delta _0)} \right)$, for some $\delta _0>0 $, we can then conclude the existence of a unique element $f _0  \in \Omega(R,W) $ such that $\Psi (f_0)= f_0 $. Now, as we already saw in the proof of part {\bf (i)}, the condition \eqref{condition for lfx in differentiable} implies that $L_{F _x}< 1 $ and hence Proposition \ref{continuous synchronization} guarantees that $f_{(\phi, \omega,F)}:M \longrightarrow W $  is the unique continuous map that satisfies the identity $\Psi (f_{(\phi, \omega,F)})= f_{(\phi, \omega,F)} $. Since the element $f _0  \in \Omega(R,W) $ will be obviously continuous, we shall immediately be able to conclude that $f_{(\phi, \omega,F)}= f _0 $ and hence that $f_{(\phi, \omega,F)} \in \Omega(R,W) \subset C ^1(M, \mathbb{R}^N) $, necessarily.

We hence conclude the proof by showing the existence of a constant $\delta_0>0 $ for which $\Psi: \Omega(R,W) \longrightarrow \Omega(R,W) $ is a contraction in the norm $\left\|\cdot \right\|_{C ^1(\delta_0)} $. Let $f _1, f _2 \in \Omega(R,W) $, and $\delta>0 $  arbitrary. By definition,
\begin{multline}
\label{psi with norm delta}
\left\|\Psi(f  _1)-\Psi(f  _2)\right\|_{C ^1(\delta)}= \left\|\Psi(f  _1)-\Psi(f  _2)\right\|_{\infty}\\
+ \delta \left\|D(\Psi(f  _1))-D(\Psi(f  _2))\right\|_{\infty}.
\end{multline}
Since in \eqref{f contraction with lfx} we showed that $F:V \times \omega(M) \longrightarrow V$ is a contraction in the state variable with constant $0<L_{F _x}<1$ then so is $F:W \times \omega(M) \longrightarrow W$. We can hence conclude using \eqref{phi contraction proof} that 
\begin{equation}
\label{psi contraction with lfx}
\!\!\!\!\!\!\!\!\!\left\|\Psi(f_1)-\Psi(f _2)\right\|_{\infty}<L_{F _x} \left\|f_1-f _2\right\|_{\infty}.
\end{equation}
Let now $m \in M $  and $\mathbf{v} \in T _mM  $ arbitrary. By the chain rule, the Mean Value Theorem, the convexity of $W$, and the definitions in \eqref{Ls for F}:
\scriptsize
\begin{multline*}
\left\|D \Psi(f_1)(m) \cdot \mathbf{v}-D \Psi(f_2)(m) \cdot \mathbf{v}\right\|\\
= \|  DF(f_1(\phi^{-1}(m)), \omega (m))\left( D f _1(\phi^{-1} (m))\cdot T _m \phi^{-1} \cdot \mathbf{v}, D\omega (m) \cdot \mathbf{v} \right)\\
-DF(f_2(\phi^{-1}(m)), \omega (m))\left( D f _2(\phi^{-1} (m))\cdot T _m \phi^{-1} \cdot \mathbf{v}, D\omega (m) \cdot \mathbf{v} \right)\|\\
\leq  \|  DF(f_1(\phi^{-1}(m)), \omega (m))\left( D f _1(\phi^{-1} (m))\cdot T _m \phi^{-1} \cdot \mathbf{v}, D\omega (m) \cdot \mathbf{v} \right)\\
-DF(f_1(\phi^{-1}(m)), \omega (m))\left( D f _2(\phi^{-1} (m))\cdot T _m \phi^{-1} \cdot \mathbf{v}, D\omega (m) \cdot \mathbf{v} \right)\|\\
+ \|  DF(f_1(\phi^{-1}(m)), \omega (m))\left( D f _2(\phi^{-1} (m))\cdot T _m \phi^{-1} \cdot \mathbf{v}, D\omega (m) \cdot \mathbf{v} \right)\\
-DF(f_2(\phi^{-1}(m)), \omega (m))\left( D f _2(\phi^{-1} (m))\cdot T _m \phi^{-1} \cdot \mathbf{v}, D\omega (m) \cdot \mathbf{v} \right)\|\\
\leq L_{F _x} \left\| \left(D f _1(\phi^{-1} (m))- D f _2(\phi^{-1} (m))\right)\cdot T _m \phi^{-1} \cdot \mathbf{v}\right\|\\
+ \|  (D _xF(f_1(\phi^{-1}(m)), \omega (m))\\
-D_xF(f_2(\phi^{-1}(m)), \omega (m))) \cdot \left(D f _2(\phi^{-1} (m))\cdot T _m \phi^{-1} \cdot \mathbf{v}\right)\| \\
+\left\|(D _zF(f_1(\phi^{-1}(m)), \omega (m))-D_zF(f_2(\phi^{-1}(m)), \omega (m))) \cdot D\omega (m) \cdot \mathbf{v}\right\|\\
\leq L_{F _x}\left\|Df _1- Df _2\right\|_{\infty}\left\|T \phi ^{-1}\right\|_{\infty} \left\|\mathbf{v}\right\|\\
+L_{F_{xx}}\left\|f _1- f _2\right\|_{\infty} \left\|D f _2\right\|_{\infty}\left\|T \phi ^{-1}\right\|_{\infty} \left\|\mathbf{v}\right\|\\
+L_{F_{x z}}\left\|f _1- f _2\right\|_{\infty} \left\|D\omega\right\|_{\infty} \left\|\mathbf{v}\right\|.
\end{multline*}
\normalsize
If we now use that, by hypothesis, $f _2 \in \Omega(R,W) $, the last inequality implies that:
\begin{multline*}
\!\!\!\!\!\!\!\left\|D \Psi(f_1)(m) \cdot \mathbf{v}-D \Psi(f_2)(m) \cdot \mathbf{v}\right\|\\
\leq (L_{F _x}\left\|Df _1- Df _2\right\|_{\infty}\left\|T \phi ^{-1}\right\|_{\infty} \\
+L_{F_{xx}}\left\|f _1- f _2\right\|_{\infty}R\left\|T \phi ^{-1}\right\|_{\infty} \\
+L_{F_{x z}}\left\|f _1- f _2\right\|_{\infty} \left\|D\omega\right\|_{\infty} ) \left\|\mathbf{v}\right\|,
\end{multline*}
which ensures that
\begin{multline*}
\!\!\!\!\left\|D \Psi(f_1)(m) \cdot \mathbf{v}-D \Psi(f_2)(m) \cdot \mathbf{v}\right\|\\
\leq(L_{F_{xx}}\left\|T \phi ^{-1}\right\|_{\infty}R+L_{F_{x z}}  \left\|D\omega\right\|_{\infty})\left\|f _1- f _2\right\|_{\infty}\\
 +L_{F _x}\left\|T \phi ^{-1}\right\|_{\infty}\left\|Df _1- Df _2\right\|_{\infty}.
\end{multline*}
Together with \eqref{definition c1d} and \eqref{psi contraction with lfx}, this inequality implies that:
\begin{multline}
\label{psic1delta}
\left\|\Psi(f  _1)-\Psi(f  _2)\right\|_{C ^1(\delta)}\\
\leq  ( L_{F _x}+ \delta \left( L_{F_{xx}}\left\|T \phi ^{-1}\right\|_{\infty}R+L_{F_{x z}}  \left\|D\omega\right\|_{\infty}  \right)  )\left\|f _1- f _2\right\|_{\infty}\\
+ \delta L_{F _x}\left\|T \phi ^{-1}\right\|_{\infty}\left\|Df _1- Df _2\right\|_{\infty}.
\end{multline}
Choose now $\delta_0>0 $ small enough so that 
\begin{equation}
\label{almost there dps}
 L_{F _x}+ \delta _0 \left( L_{F_{xx}}\left\|T \phi ^{-1}\right\|_{\infty}R+L_{F_{x z}}  \left\|D\omega\right\|_{\infty}  \right)<1,
\end{equation}
which is always possible due to the hypothesis $L_{F _x}<1 $. Additionally, the hypothesis $L_{F _x} \left\|T \phi^{-1}\right\|_{\infty}<1 $ together with \eqref{almost there dps} and \eqref{psic1delta} imply that 
\begin{equation*}
\left\|\Psi(f  _1)-\Psi(f  _2)\right\|_{C ^1(\delta_0)}\leq c _0
\left\|f  _1- f  _2\right\|_{C ^1(\delta_0)},
\end{equation*}
for any $0< c  _0<1  $ such that 
\begin{multline*}
c _0<\min \{L_{F _x} \left\|T \phi^{-1}\right\|_{\infty},  \\
L_{F _x}+ \delta_0 \left( L_{F_{xx}}\left\|T \phi ^{-1}\right\|_{\infty}R+L_{F_{x z}}  \left\|D\omega\right\|_{\infty}  \right)\},
\end{multline*}
which guarantees the contractivity of $\Psi: \Omega(R,W) \longrightarrow \Omega(R,W) $, as required. \quad $\blacksquare$

\medskip


\begin{acknowledgments}
We thank Edward Ott for pointing out the relation of our research in reservoir computing with the differentiability problem. His observation is at the origin of this contribution. We also thank Herbert Jaeger for his comment about global versus local contractivity as well as the suggestions of the referees that have significantly improved the paper.

AH is supported by a scholarship from the EPSRC Centre for Doctoral Training in Statistical Applied Mathematics at Bath (SAMBa), project EP/L015684/1. JPO acknowledges partial financial support  coming from the Research Commission of the Universit\"at Sankt Gallen, the Swiss National Science Foundation (grant number 200021\_175801/1), and the French ANR ``BIPHOPROC" project (ANR-14-OHRI-0002-02).  The authors thank the hospitality and the generosity of the FIM at ETH Zurich and the Division of Mathematical Sciences of the Nanyang Technological University, Singapore, where a significant portion of the results in this paper were obtained. 
\end{acknowledgments}


\begin{thebibliography}{53}%
\makeatletter
\providecommand \@ifxundefined [1]{%
 \@ifx{#1\undefined}
}%
\providecommand \@ifnum [1]{%
 \ifnum #1\expandafter \@firstoftwo
 \else \expandafter \@secondoftwo
 \fi
}%
\providecommand \@ifx [1]{%
 \ifx #1\expandafter \@firstoftwo
 \else \expandafter \@secondoftwo
 \fi
}%
\providecommand \natexlab [1]{#1}%
\providecommand \enquote  [1]{``#1''}%
\providecommand \bibnamefont  [1]{#1}%
\providecommand \bibfnamefont [1]{#1}%
\providecommand \citenamefont [1]{#1}%
\providecommand \href@noop [0]{\@secondoftwo}%
\providecommand \href [0]{\begingroup \@sanitize@url \@href}%
\providecommand \@href[1]{\@@startlink{#1}\@@href}%
\providecommand \@@href[1]{\endgroup#1\@@endlink}%
\providecommand \@sanitize@url [0]{\catcode `\\12\catcode `\$12\catcode
  `\&12\catcode `\#12\catcode `\^12\catcode `\_12\catcode `\%12\relax}%
\providecommand \@@startlink[1]{}%
\providecommand \@@endlink[0]{}%
\providecommand \url  [0]{\begingroup\@sanitize@url \@url }%
\providecommand \@url [1]{\endgroup\@href {#1}{\urlprefix }}%
\providecommand \urlprefix  [0]{URL }%
\providecommand \Eprint [0]{\href }%
\providecommand \doibase [0]{https://doi.org/}%
\providecommand \selectlanguage [0]{\@gobble}%
\providecommand \bibinfo  [0]{\@secondoftwo}%
\providecommand \bibfield  [0]{\@secondoftwo}%
\providecommand \translation [1]{[#1]}%
\providecommand \BibitemOpen [0]{}%
\providecommand \bibitemStop [0]{}%
\providecommand \bibitemNoStop [0]{.\EOS\space}%
\providecommand \EOS [0]{\spacefactor3000\relax}%
\providecommand \BibitemShut  [1]{\csname bibitem#1\endcsname}%
\let\auto@bib@innerbib\@empty
\bibitem [{\citenamefont {Pecora}\ \emph {et~al.}(1997)\citenamefont {Pecora},
  \citenamefont {Carroll}, \citenamefont {Johnson}, \citenamefont {Mar},\ and\
  \citenamefont {Heagy}}]{pecora:synch}%
  \BibitemOpen
  \bibfield  {author} {\bibinfo {author} {\bibfnamefont {L.~M.}\ \bibnamefont
  {Pecora}}, \bibinfo {author} {\bibfnamefont {T.~L.}\ \bibnamefont {Carroll}},
  \bibinfo {author} {\bibfnamefont {G.~A.}\ \bibnamefont {Johnson}}, \bibinfo
  {author} {\bibfnamefont {D.~J.}\ \bibnamefont {Mar}},\ and\ \bibinfo {author}
  {\bibfnamefont {J.~F.}\ \bibnamefont {Heagy}},\ }\bibfield  {title} {\bibinfo
  {title} {{Fundamentals of synchronization in chaotic systems, concepts, and
  applications}},\ }\href {https://doi.org/10.1063/1.166278} {\bibfield
  {journal} {\bibinfo  {journal} {Chaos}\ }\textbf {\bibinfo {volume} {7}},\
  \bibinfo {pages} {520} (\bibinfo {year} {1997})}\BibitemShut {NoStop}%
\bibitem [{\citenamefont {Ott}(2002)}]{ott2002chaos}%
  \BibitemOpen
  \bibfield  {author} {\bibinfo {author} {\bibfnamefont {E.}~\bibnamefont
  {Ott}},\ }\href@noop {} {\emph {\bibinfo {title} {{Chaos in Dynamical
  Systems}}}},\ \bibinfo {edition} {2nd}\ ed.\ (\bibinfo  {publisher}
  {Cambridge University Press},\ \bibinfo {year} {2002})\BibitemShut {NoStop}%
\bibitem [{\citenamefont {Boccaletti}\ \emph {et~al.}(2002)\citenamefont
  {Boccaletti}, \citenamefont {Kurths}, \citenamefont {Osipov}, \citenamefont
  {Valladares},\ and\ \citenamefont {Zhou}}]{boccaletti:reports:2002}%
  \BibitemOpen
  \bibfield  {author} {\bibinfo {author} {\bibfnamefont {S.}~\bibnamefont
  {Boccaletti}}, \bibinfo {author} {\bibfnamefont {J.}~\bibnamefont {Kurths}},
  \bibinfo {author} {\bibfnamefont {G.}~\bibnamefont {Osipov}}, \bibinfo
  {author} {\bibfnamefont {D.~L.}\ \bibnamefont {Valladares}},\ and\ \bibinfo
  {author} {\bibfnamefont {C.~S.}\ \bibnamefont {Zhou}},\ }\bibfield  {title}
  {\bibinfo {title} {{The synchronization of chaotic systems}},\ }\href@noop {}
  {\bibfield  {journal} {\bibinfo  {journal} {Physics Reports}\ }\textbf
  {\bibinfo {volume} {366}},\ \bibinfo {pages} {1} (\bibinfo {year}
  {2002})}\BibitemShut {NoStop}%
\bibitem [{\citenamefont {Eroglu}\ \emph {et~al.}(2017)\citenamefont {Eroglu},
  \citenamefont {Lamb},\ and\ \citenamefont
  {Pereira}}]{eroglu2017synchronisation}%
  \BibitemOpen
  \bibfield  {author} {\bibinfo {author} {\bibfnamefont {D.}~\bibnamefont
  {Eroglu}}, \bibinfo {author} {\bibfnamefont {J.~S.~W.}\ \bibnamefont
  {Lamb}},\ and\ \bibinfo {author} {\bibfnamefont {T.}~\bibnamefont
  {Pereira}},\ }\bibfield  {title} {\bibinfo {title} {{Synchronisation of chaos
  and its applications}},\ }\href@noop {} {\bibfield  {journal} {\bibinfo
  {journal} {Contemporary Physics}\ }\textbf {\bibinfo {volume} {58}},\
  \bibinfo {pages} {207} (\bibinfo {year} {2017})}\BibitemShut {NoStop}%
\bibitem [{\citenamefont {Alvarez}\ \emph {et~al.}(2005)\citenamefont
  {Alvarez}, \citenamefont {Li}, \citenamefont {Montoya}, \citenamefont
  {Pastor},\ and\ \citenamefont {Romera}}]{alvarez2005breaking}%
  \BibitemOpen
  \bibfield  {author} {\bibinfo {author} {\bibfnamefont {G.}~\bibnamefont
  {Alvarez}}, \bibinfo {author} {\bibfnamefont {S.}~\bibnamefont {Li}},
  \bibinfo {author} {\bibfnamefont {F.}~\bibnamefont {Montoya}}, \bibinfo
  {author} {\bibfnamefont {G.}~\bibnamefont {Pastor}},\ and\ \bibinfo {author}
  {\bibfnamefont {M.}~\bibnamefont {Romera}},\ }\bibfield  {title} {\bibinfo
  {title} {{Breaking projective chaos synchronization secure communication
  using filtering and generalized synchronization}},\ }\href@noop {} {\bibfield
   {journal} {\bibinfo  {journal} {Chaos, Solitons {\&} Fractals}\ }\textbf
  {\bibinfo {volume} {24}},\ \bibinfo {pages} {775} (\bibinfo {year}
  {2005})}\BibitemShut {NoStop}%
\bibitem [{\citenamefont {Moskalenko}\ \emph {et~al.}(2010)\citenamefont
  {Moskalenko}, \citenamefont {Koronovskii},\ and\ \citenamefont
  {Hramov}}]{moskalenko2010generalized}%
  \BibitemOpen
  \bibfield  {author} {\bibinfo {author} {\bibfnamefont {O.~I.}\ \bibnamefont
  {Moskalenko}}, \bibinfo {author} {\bibfnamefont {A.~A.}\ \bibnamefont
  {Koronovskii}},\ and\ \bibinfo {author} {\bibfnamefont {A.~E.}\ \bibnamefont
  {Hramov}},\ }\bibfield  {title} {\bibinfo {title} {{Generalized
  synchronization of chaos for secure communication: Remarkable stability to
  noise}},\ }\href@noop {} {\bibfield  {journal} {\bibinfo  {journal} {Physics
  Letters A}\ }\textbf {\bibinfo {volume} {374}},\ \bibinfo {pages} {2925}
  (\bibinfo {year} {2010})}\BibitemShut {NoStop}%
\bibitem [{\citenamefont {Jiang-Feng}\ \emph {et~al.}(2004)\citenamefont
  {Jiang-Feng}, \citenamefont {Le-Quan},\ and\ \citenamefont
  {Guan-Rong}}]{jiang2004chaotic}%
  \BibitemOpen
  \bibfield  {author} {\bibinfo {author} {\bibfnamefont {X.}~\bibnamefont
  {Jiang-Feng}}, \bibinfo {author} {\bibfnamefont {M.}~\bibnamefont
  {Le-Quan}},\ and\ \bibinfo {author} {\bibfnamefont {C.}~\bibnamefont
  {Guan-Rong}},\ }\bibfield  {title} {\bibinfo {title} {{A chaotic
  communication scheme based on generalized synchronization and hash
  functions}},\ }\href@noop {} {\bibfield  {journal} {\bibinfo  {journal}
  {Chinese Physics Letters}\ }\textbf {\bibinfo {volume} {21}},\ \bibinfo
  {pages} {1445} (\bibinfo {year} {2004})}\BibitemShut {NoStop}%
\bibitem [{\citenamefont {Stam}\ \emph {et~al.}(2002)\citenamefont {Stam},
  \citenamefont {van Walsum}, \citenamefont {Pijnenburg}, \citenamefont
  {Berendse}, \citenamefont {de~Munck}, \citenamefont {Scheltens},\ and\
  \citenamefont {van Dijk}}]{stam2002generalized}%
  \BibitemOpen
  \bibfield  {author} {\bibinfo {author} {\bibfnamefont {C.~J.}\ \bibnamefont
  {Stam}}, \bibinfo {author} {\bibfnamefont {A.~M. v.~C.}\ \bibnamefont {van
  Walsum}}, \bibinfo {author} {\bibfnamefont {Y.~A.~L.}\ \bibnamefont
  {Pijnenburg}}, \bibinfo {author} {\bibfnamefont {H.~W.}\ \bibnamefont
  {Berendse}}, \bibinfo {author} {\bibfnamefont {J.~C.}\ \bibnamefont
  {de~Munck}}, \bibinfo {author} {\bibfnamefont {P.}~\bibnamefont
  {Scheltens}},\ and\ \bibinfo {author} {\bibfnamefont {B.~W.}\ \bibnamefont
  {van Dijk}},\ }\bibfield  {title} {\bibinfo {title} {{Generalized
  synchronization of MEG recordings in Alzheimer's disease: evidence for
  involvement of the gamma band}},\ }\href@noop {} {\bibfield  {journal}
  {\bibinfo  {journal} {Journal of Clinical Neurophysiology}\ }\textbf
  {\bibinfo {volume} {19}},\ \bibinfo {pages} {562} (\bibinfo {year}
  {2002})}\BibitemShut {NoStop}%
\bibitem [{\citenamefont {Bartolomei}\ \emph {et~al.}(2006)\citenamefont
  {Bartolomei}, \citenamefont {Bosma}, \citenamefont {Klein}, \citenamefont
  {Baayen}, \citenamefont {Reijneveld}, \citenamefont {Postma}, \citenamefont
  {Heimans}, \citenamefont {van Dijk}, \citenamefont {de~Munck}, \citenamefont
  {de~Jongh}, \citenamefont {Cover},\ and\ \citenamefont
  {Stam}}]{bartolomei2006disturbed}%
  \BibitemOpen
  \bibfield  {author} {\bibinfo {author} {\bibfnamefont {F.}~\bibnamefont
  {Bartolomei}}, \bibinfo {author} {\bibfnamefont {I.}~\bibnamefont {Bosma}},
  \bibinfo {author} {\bibfnamefont {M.}~\bibnamefont {Klein}}, \bibinfo
  {author} {\bibfnamefont {J.~C.}\ \bibnamefont {Baayen}}, \bibinfo {author}
  {\bibfnamefont {J.~C.}\ \bibnamefont {Reijneveld}}, \bibinfo {author}
  {\bibfnamefont {T.~J.}\ \bibnamefont {Postma}}, \bibinfo {author}
  {\bibfnamefont {J.~J.}\ \bibnamefont {Heimans}}, \bibinfo {author}
  {\bibfnamefont {B.~W.}\ \bibnamefont {van Dijk}}, \bibinfo {author}
  {\bibfnamefont {J.~C.}\ \bibnamefont {de~Munck}}, \bibinfo {author}
  {\bibfnamefont {A.}~\bibnamefont {de~Jongh}}, \bibinfo {author}
  {\bibfnamefont {K.~S.}\ \bibnamefont {Cover}},\ and\ \bibinfo {author}
  {\bibfnamefont {C.~J.}\ \bibnamefont {Stam}},\ }\bibfield  {title} {\bibinfo
  {title} {{Disturbed functional connectivity in brain tumour patients:
  evaluation by graph analysis of synchronization matrices}},\ }\href@noop {}
  {\bibfield  {journal} {\bibinfo  {journal} {Clinical Neurophysiology}\
  }\textbf {\bibinfo {volume} {117}},\ \bibinfo {pages} {2039} (\bibinfo {year}
  {2006})}\BibitemShut {NoStop}%
\bibitem [{\citenamefont {Rulkov}\ \emph {et~al.}(1995)\citenamefont {Rulkov},
  \citenamefont {Sushchik}, \citenamefont {Tsimring},\ and\ \citenamefont
  {Abarbanel}}]{rulkov1995generalized}%
  \BibitemOpen
  \bibfield  {author} {\bibinfo {author} {\bibfnamefont {N.~F.}\ \bibnamefont
  {Rulkov}}, \bibinfo {author} {\bibfnamefont {M.~M.}\ \bibnamefont
  {Sushchik}}, \bibinfo {author} {\bibfnamefont {L.~S.}\ \bibnamefont
  {Tsimring}},\ and\ \bibinfo {author} {\bibfnamefont {H.~D.~I.}\ \bibnamefont
  {Abarbanel}},\ }\bibfield  {title} {\bibinfo {title} {{Generalized
  synchronization of chaos in directionally coupled chaotic systems}},\
  }\href@noop {} {\bibfield  {journal} {\bibinfo  {journal} {Physical Review
  E}\ }\textbf {\bibinfo {volume} {51}},\ \bibinfo {pages} {980} (\bibinfo
  {year} {1995})}\BibitemShut {NoStop}%
\bibitem [{\citenamefont {Kocarev}\ and\ \citenamefont
  {Parlitz}(1995)}]{kocarev1995general}%
  \BibitemOpen
  \bibfield  {author} {\bibinfo {author} {\bibfnamefont {L.}~\bibnamefont
  {Kocarev}}\ and\ \bibinfo {author} {\bibfnamefont {U.}~\bibnamefont
  {Parlitz}},\ }\bibfield  {title} {\bibinfo {title} {{General approach for
  chaotic synchronization with applications to communication}},\ }\href@noop {}
  {\bibfield  {journal} {\bibinfo  {journal} {Physical Review Letters}\
  }\textbf {\bibinfo {volume} {74}},\ \bibinfo {pages} {5028} (\bibinfo {year}
  {1995})}\BibitemShut {NoStop}%
\bibitem [{\citenamefont {Pyragas}(1996)}]{pyragas:1996}%
  \BibitemOpen
  \bibfield  {author} {\bibinfo {author} {\bibfnamefont {K.}~\bibnamefont
  {Pyragas}},\ }\bibfield  {title} {\bibinfo {title} {{Weak and strong
  synchronization of chaos}},\ }\href
  {https://doi.org/10.1103/PhysRevE.54.R4508} {\bibfield  {journal} {\bibinfo
  {journal} {Physical Review E - Statistical Physics, Plasmas, Fluids, and
  Related Interdisciplinary Topics}\ }\textbf {\bibinfo {volume} {54}},\
  \bibinfo {pages} {4508} (\bibinfo {year} {1996})}\BibitemShut {NoStop}%
\bibitem [{\citenamefont {Hunt}\ \emph {et~al.}(1997)\citenamefont {Hunt},
  \citenamefont {Ott},\ and\ \citenamefont {Yorke}}]{hunt:ott:1997}%
  \BibitemOpen
  \bibfield  {author} {\bibinfo {author} {\bibfnamefont {B.~R.}\ \bibnamefont
  {Hunt}}, \bibinfo {author} {\bibfnamefont {E.}~\bibnamefont {Ott}},\ and\
  \bibinfo {author} {\bibfnamefont {J.~A.}\ \bibnamefont {Yorke}},\ }\bibfield
  {title} {\bibinfo {title} {{Differentiable generalized synchronization of
  chaos}},\ }\href {https://doi.org/10.1103/PhysRevE.55.4029} {\bibfield
  {journal} {\bibinfo  {journal} {Physical Review E}\ }\textbf {\bibinfo
  {volume} {55}},\ \bibinfo {pages} {4029} (\bibinfo {year}
  {1997})}\BibitemShut {NoStop}%
\bibitem [{\citenamefont {Luko{\v{s}}evi{\v{c}}ius}\ and\ \citenamefont
  {Jaeger}(2009)}]{lukosevicius}%
  \BibitemOpen
  \bibfield  {author} {\bibinfo {author} {\bibfnamefont {M.}~\bibnamefont
  {Luko{\v{s}}evi{\v{c}}ius}}\ and\ \bibinfo {author} {\bibfnamefont
  {H.}~\bibnamefont {Jaeger}},\ }\bibfield  {title} {\bibinfo {title}
  {{Reservoir computing approaches to recurrent neural network training}},\
  }\href@noop {} {\bibfield  {journal} {\bibinfo  {journal} {Computer Science
  Review}\ }\textbf {\bibinfo {volume} {3}},\ \bibinfo {pages} {127} (\bibinfo
  {year} {2009})}\BibitemShut {NoStop}%
\bibitem [{\citenamefont {Tanaka}\ \emph {et~al.}(2019)\citenamefont {Tanaka},
  \citenamefont {Yamane}, \citenamefont {H{\'{e}}roux}, \citenamefont {Nakane},
  \citenamefont {Kanazawa}, \citenamefont {Takeda}, \citenamefont {Numata},
  \citenamefont {Nakano},\ and\ \citenamefont {Hirose}}]{tanaka:review}%
  \BibitemOpen
  \bibfield  {author} {\bibinfo {author} {\bibfnamefont {G.}~\bibnamefont
  {Tanaka}}, \bibinfo {author} {\bibfnamefont {T.}~\bibnamefont {Yamane}},
  \bibinfo {author} {\bibfnamefont {J.~B.}\ \bibnamefont {H{\'{e}}roux}},
  \bibinfo {author} {\bibfnamefont {R.}~\bibnamefont {Nakane}}, \bibinfo
  {author} {\bibfnamefont {N.}~\bibnamefont {Kanazawa}}, \bibinfo {author}
  {\bibfnamefont {S.}~\bibnamefont {Takeda}}, \bibinfo {author} {\bibfnamefont
  {H.}~\bibnamefont {Numata}}, \bibinfo {author} {\bibfnamefont
  {D.}~\bibnamefont {Nakano}},\ and\ \bibinfo {author} {\bibfnamefont
  {A.}~\bibnamefont {Hirose}},\ }\bibfield  {title} {\bibinfo {title} {{Recent
  advances in physical reservoir computing: A review}},\ }\href
  {https://doi.org/10.1016/j.neunet.2019.03.005} {\bibfield  {journal}
  {\bibinfo  {journal} {Neural Networks}\ }\textbf {\bibinfo {volume} {115}},\
  \bibinfo {pages} {100} (\bibinfo {year} {2019})}\BibitemShut {NoStop}%
\bibitem [{\citenamefont {Jaeger}\ and\ \citenamefont {Haas}(2004)}]{Jaeger04}%
  \BibitemOpen
  \bibfield  {author} {\bibinfo {author} {\bibfnamefont {H.}~\bibnamefont
  {Jaeger}}\ and\ \bibinfo {author} {\bibfnamefont {H.}~\bibnamefont {Haas}},\
  }\bibfield  {title} {\bibinfo {title} {{Harnessing Nonlinearity: Predicting
  Chaotic Systems and Saving Energy in Wireless Communication}},\ }\href
  {https://doi.org/10.1126/science.1091277} {\bibfield  {journal} {\bibinfo
  {journal} {Science}\ }\textbf {\bibinfo {volume} {304}},\ \bibinfo {pages}
  {78} (\bibinfo {year} {2004})}\BibitemShut {NoStop}%
\bibitem [{\citenamefont {Pathak}\ \emph {et~al.}(2017)\citenamefont {Pathak},
  \citenamefont {Lu}, \citenamefont {Hunt}, \citenamefont {Girvan},\ and\
  \citenamefont {Ott}}]{pathak:chaos}%
  \BibitemOpen
  \bibfield  {author} {\bibinfo {author} {\bibfnamefont {J.}~\bibnamefont
  {Pathak}}, \bibinfo {author} {\bibfnamefont {Z.}~\bibnamefont {Lu}}, \bibinfo
  {author} {\bibfnamefont {B.~R.}\ \bibnamefont {Hunt}}, \bibinfo {author}
  {\bibfnamefont {M.}~\bibnamefont {Girvan}},\ and\ \bibinfo {author}
  {\bibfnamefont {E.}~\bibnamefont {Ott}},\ }\bibfield  {title} {\bibinfo
  {title} {{Using machine learning to replicate chaotic attractors and
  calculate Lyapunov exponents from data}},\ }\bibfield  {journal} {\bibinfo
  {journal} {Chaos}\ }\textbf {\bibinfo {volume} {27}},\ \href
  {https://doi.org/10.1063/1.5010300} {10.1063/1.5010300} (\bibinfo {year}
  {2017}),\ \Eprint {https://arxiv.org/abs/1710.07313} {arXiv:1710.07313}
  \BibitemShut {NoStop}%
\bibitem [{\citenamefont {Pathak}\ \emph {et~al.}(2018)\citenamefont {Pathak},
  \citenamefont {Hunt}, \citenamefont {Girvan}, \citenamefont {Lu},\ and\
  \citenamefont {Ott}}]{Pathak:PRL}%
  \BibitemOpen
  \bibfield  {author} {\bibinfo {author} {\bibfnamefont {J.}~\bibnamefont
  {Pathak}}, \bibinfo {author} {\bibfnamefont {B.}~\bibnamefont {Hunt}},
  \bibinfo {author} {\bibfnamefont {M.}~\bibnamefont {Girvan}}, \bibinfo
  {author} {\bibfnamefont {Z.}~\bibnamefont {Lu}},\ and\ \bibinfo {author}
  {\bibfnamefont {E.}~\bibnamefont {Ott}},\ }\bibfield  {title} {\bibinfo
  {title} {{Model-Free Prediction of Large Spatiotemporally Chaotic Systems
  from Data: A Reservoir Computing Approach}},\ }\href
  {https://doi.org/10.1103/PhysRevLett.120.024102} {\bibfield  {journal}
  {\bibinfo  {journal} {Physical Review Letters}\ }\textbf {\bibinfo {volume}
  {120}},\ \bibinfo {pages} {24102} (\bibinfo {year} {2018})}\BibitemShut
  {NoStop}%
\bibitem [{\citenamefont {Lu}\ \emph {et~al.}(2018)\citenamefont {Lu},
  \citenamefont {Hunt},\ and\ \citenamefont {Ott}}]{Ott2018}%
  \BibitemOpen
  \bibfield  {author} {\bibinfo {author} {\bibfnamefont {Z.}~\bibnamefont
  {Lu}}, \bibinfo {author} {\bibfnamefont {B.~R.}\ \bibnamefont {Hunt}},\ and\
  \bibinfo {author} {\bibfnamefont {E.}~\bibnamefont {Ott}},\ }\bibfield
  {title} {\bibinfo {title} {{Attractor reconstruction by machine learning}},\
  }\bibfield  {journal} {\bibinfo  {journal} {Chaos}\ }\textbf {\bibinfo
  {volume} {28}},\ \href {https://doi.org/10.1063/1.5039508}
  {10.1063/1.5039508} (\bibinfo {year} {2018}),\ \Eprint
  {https://arxiv.org/abs/1805.03362} {arXiv:1805.03362} \BibitemShut {NoStop}%
\bibitem [{\citenamefont {Carroll}(2018)}]{carroll2018using}%
  \BibitemOpen
  \bibfield  {author} {\bibinfo {author} {\bibfnamefont {T.~L.}\ \bibnamefont
  {Carroll}},\ }\bibfield  {title} {\bibinfo {title} {{Using reservoir
  computers to distinguish chaotic signals}},\ }\href@noop {} {\bibfield
  {journal} {\bibinfo  {journal} {Physical Review E}\ }\textbf {\bibinfo
  {volume} {98}},\ \bibinfo {pages} {52209} (\bibinfo {year}
  {2018})}\BibitemShut {NoStop}%
\bibitem [{\citenamefont {Ib{\'{a}}{\~{n}}ez-Soria}\ \emph
  {et~al.}(2018)\citenamefont {Ib{\'{a}}{\~{n}}ez-Soria}, \citenamefont
  {Garc{\'{i}}a-Ojalvo}, \citenamefont {Soria-Frisch},\ and\ \citenamefont
  {Ruffini}}]{ibanez2018detection}%
  \BibitemOpen
  \bibfield  {author} {\bibinfo {author} {\bibfnamefont {D.}~\bibnamefont
  {Ib{\'{a}}{\~{n}}ez-Soria}}, \bibinfo {author} {\bibfnamefont
  {J.}~\bibnamefont {Garc{\'{i}}a-Ojalvo}}, \bibinfo {author} {\bibfnamefont
  {A.}~\bibnamefont {Soria-Frisch}},\ and\ \bibinfo {author} {\bibfnamefont
  {G.}~\bibnamefont {Ruffini}},\ }\bibfield  {title} {\bibinfo {title}
  {{Detection of generalized synchronization using echo state networks}},\
  }\href@noop {} {\bibfield  {journal} {\bibinfo  {journal} {Chaos: An
  Interdisciplinary Journal of Nonlinear Science}\ }\textbf {\bibinfo {volume}
  {28}},\ \bibinfo {pages} {33118} (\bibinfo {year} {2018})}\BibitemShut
  {NoStop}%
\bibitem [{\citenamefont {Weng}\ \emph {et~al.}(2019)\citenamefont {Weng},
  \citenamefont {Yang}, \citenamefont {Gu}, \citenamefont {Zhang},\ and\
  \citenamefont {Small}}]{Weng:2019}%
  \BibitemOpen
  \bibfield  {author} {\bibinfo {author} {\bibfnamefont {T.}~\bibnamefont
  {Weng}}, \bibinfo {author} {\bibfnamefont {H.}~\bibnamefont {Yang}}, \bibinfo
  {author} {\bibfnamefont {C.}~\bibnamefont {Gu}}, \bibinfo {author}
  {\bibfnamefont {J.}~\bibnamefont {Zhang}},\ and\ \bibinfo {author}
  {\bibfnamefont {M.}~\bibnamefont {Small}},\ }\bibfield  {title} {\bibinfo
  {title} {{Synchronization of chaotic systems and their machine-learning
  models}},\ }\href {https://doi.org/10.1103/PhysRevE.99.042203} {\bibfield
  {journal} {\bibinfo  {journal} {Physical Review E}\ }\textbf {\bibinfo
  {volume} {99}},\ \bibinfo {pages} {1} (\bibinfo {year} {2019})}\BibitemShut
  {NoStop}%
\bibitem [{\citenamefont {Lymburn}\ \emph {et~al.}(2019)\citenamefont
  {Lymburn}, \citenamefont {Walker}, \citenamefont {Small},\ and\ \citenamefont
  {J{\"{u}}ngling}}]{Lymburn:2019}%
  \BibitemOpen
  \bibfield  {author} {\bibinfo {author} {\bibfnamefont {T.}~\bibnamefont
  {Lymburn}}, \bibinfo {author} {\bibfnamefont {D.~M.}\ \bibnamefont {Walker}},
  \bibinfo {author} {\bibfnamefont {M.}~\bibnamefont {Small}},\ and\ \bibinfo
  {author} {\bibfnamefont {T.}~\bibnamefont {J{\"{u}}ngling}},\ }\bibfield
  {title} {\bibinfo {title} {{The reservoir's perspective on generalized
  synchronization}},\ }\bibfield  {journal} {\bibinfo  {journal} {Chaos}\
  }\textbf {\bibinfo {volume} {29}},\ \href {https://doi.org/10.1063/1.5120733}
  {10.1063/1.5120733} (\bibinfo {year} {2019})\BibitemShut {NoStop}%
\bibitem [{\citenamefont {Takens}(1981)}]{takensembedding}%
  \BibitemOpen
  \bibfield  {author} {\bibinfo {author} {\bibfnamefont {F.}~\bibnamefont
  {Takens}},\ }\bibfield  {title} {\bibinfo {title} {{Detecting strange
  attractors in turbulence}}\ }(\bibinfo  {publisher} {Springer Berlin
  Heidelberg},\ \bibinfo {year} {1981})\ pp.\ \bibinfo {pages}
  {366--381}\BibitemShut {NoStop}%
\bibitem [{\citenamefont {Huke}(2006)}]{huke:2006}%
  \BibitemOpen
  \bibfield  {author} {\bibinfo {author} {\bibfnamefont {J.~P.}\ \bibnamefont
  {Huke}},\ }\href@noop {} {\emph {\bibinfo {title} {{Embedding nonlinear
  dynamical systems: a guide to Takens' theorem}}}},\ \bibinfo {type} {Tech.
  Rep.}\ (\bibinfo  {institution} {Manchester Institute for Mathematical
  Sciences. The University of Manchester},\ \bibinfo {year} {2006})\BibitemShut
  {NoStop}%
\bibitem [{\citenamefont {Hart}\ \emph {et~al.}(2020)\citenamefont {Hart},
  \citenamefont {Hook},\ and\ \citenamefont {Dawes}}]{hart:ESNs}%
  \BibitemOpen
  \bibfield  {author} {\bibinfo {author} {\bibfnamefont {A.~G.}\ \bibnamefont
  {Hart}}, \bibinfo {author} {\bibfnamefont {J.~L.}\ \bibnamefont {Hook}},\
  and\ \bibinfo {author} {\bibfnamefont {J.~H.~P.}\ \bibnamefont {Dawes}},\
  }\bibfield  {title} {\bibinfo {title} {{Embedding and approximation theorems
  for echo state networks}},\ }\href {http://arxiv.org/abs/1908.05202}
  {\bibfield  {journal} {\bibinfo  {journal} {Neural Networks}\ }\textbf
  {\bibinfo {volume} {128}},\ \bibinfo {pages} {234} (\bibinfo {year}
  {2020})},\ \Eprint {https://arxiv.org/abs/1908.05202} {arXiv:1908.05202}
  \BibitemShut {NoStop}%
\bibitem [{\citenamefont {Matthews}(1992)}]{Matthews:thesis}%
  \BibitemOpen
  \bibfield  {author} {\bibinfo {author} {\bibfnamefont {M.~B.}\ \bibnamefont
  {Matthews}},\ }\emph {\bibinfo {title} {{On the Uniform Approximation of
  Nonlinear Discrete-Time Fading-Memory Systems Using Neural Network
  Models}}},\ \href {https://doi.org/10.3929/ETHZ-A-000625223} {Ph.D. thesis},\
  \bibinfo  {school} {ETH Z{\"{u}}rich} (\bibinfo {year} {1992})\BibitemShut
  {NoStop}%
\bibitem [{\citenamefont {Matthews}(1993)}]{Matthews1993}%
  \BibitemOpen
  \bibfield  {author} {\bibinfo {author} {\bibfnamefont {M.~B.}\ \bibnamefont
  {Matthews}},\ }\bibfield  {title} {\bibinfo {title} {{Approximating nonlinear
  fading-memory operators using neural network models}},\ }\href
  {https://doi.org/10.1007/BF01189878} {\bibfield  {journal} {\bibinfo
  {journal} {Circuits, Systems, and Signal Processing}\ }\textbf {\bibinfo
  {volume} {12}},\ \bibinfo {pages} {279} (\bibinfo {year} {1993})}\BibitemShut
  {NoStop}%
\bibitem [{\citenamefont {Grigoryeva}\ and\ \citenamefont
  {Ortega}(2018{\natexlab{a}})}]{RC7}%
  \BibitemOpen
  \bibfield  {author} {\bibinfo {author} {\bibfnamefont {L.}~\bibnamefont
  {Grigoryeva}}\ and\ \bibinfo {author} {\bibfnamefont {J.-P.}\ \bibnamefont
  {Ortega}},\ }\bibfield  {title} {\bibinfo {title} {{Echo state networks are
  universal}},\ }\href@noop {} {\bibfield  {journal} {\bibinfo  {journal}
  {Neural Networks}\ }\textbf {\bibinfo {volume} {108}},\ \bibinfo {pages}
  {495} (\bibinfo {year} {2018}{\natexlab{a}})},\ \Eprint
  {https://arxiv.org/abs//arxiv.org/abs/1806.00797}
  {arXiv:/arxiv.org/abs/1806.00797 [http:]} \BibitemShut {NoStop}%
\bibitem [{\citenamefont {Gonon}\ and\ \citenamefont {Ortega}(2020)}]{RC8}%
  \BibitemOpen
  \bibfield  {author} {\bibinfo {author} {\bibfnamefont {L.}~\bibnamefont
  {Gonon}}\ and\ \bibinfo {author} {\bibfnamefont {J.-P.}\ \bibnamefont
  {Ortega}},\ }\bibfield  {title} {\bibinfo {title} {{Reservoir computing
  universality with stochastic inputs}},\ }\href@noop {} {\bibfield  {journal}
  {\bibinfo  {journal} {IEEE Transactions on Neural Networks and Learning
  Systems}\ }\textbf {\bibinfo {volume} {31}},\ \bibinfo {pages} {100}
  (\bibinfo {year} {2020})},\ \Eprint {https://arxiv.org/abs/1807.02621}
  {arXiv:1807.02621} \BibitemShut {NoStop}%
\bibitem [{\citenamefont {Poggio}\ \emph {et~al.}(2017)\citenamefont {Poggio},
  \citenamefont {Mhaskar}, \citenamefont {Rosasco}, \citenamefont {Miranda},\
  and\ \citenamefont {Liao}}]{Poggio2017}%
  \BibitemOpen
  \bibfield  {author} {\bibinfo {author} {\bibfnamefont {T.}~\bibnamefont
  {Poggio}}, \bibinfo {author} {\bibfnamefont {H.}~\bibnamefont {Mhaskar}},
  \bibinfo {author} {\bibfnamefont {L.}~\bibnamefont {Rosasco}}, \bibinfo
  {author} {\bibfnamefont {B.}~\bibnamefont {Miranda}},\ and\ \bibinfo {author}
  {\bibfnamefont {Q.}~\bibnamefont {Liao}},\ }\href
  {https://doi.org/10.1007/s11633-017-1054-2} {\bibinfo {title} {{Why and when
  can deep-but not shallow-networks avoid the curse of dimensionality: A
  review}}} (\bibinfo {year} {2017}),\ \Eprint
  {https://arxiv.org/abs/1611.00740} {arXiv:1611.00740} \BibitemShut {NoStop}%
\bibitem [{\citenamefont {Mhaskar}(1996)}]{Mhaskar1996}%
  \BibitemOpen
  \bibfield  {author} {\bibinfo {author} {\bibfnamefont {N.~H.}\ \bibnamefont
  {Mhaskar}},\ }\bibfield  {title} {\bibinfo {title} {{Neural networks for
  optimal approximation of smooth and analytic functions}},\ }\href@noop {}
  {\bibfield  {journal} {\bibinfo  {journal} {Neural computation}\ }\textbf
  {\bibinfo {volume} {8}},\ \bibinfo {pages} {164} (\bibinfo {year}
  {1996})}\BibitemShut {NoStop}%
\bibitem [{\citenamefont {Inubushi}\ and\ \citenamefont
  {Goto}(2020)}]{inubushi2020transfer}%
  \BibitemOpen
  \bibfield  {author} {\bibinfo {author} {\bibfnamefont {M.}~\bibnamefont
  {Inubushi}}\ and\ \bibinfo {author} {\bibfnamefont {S.}~\bibnamefont
  {Goto}},\ }\bibfield  {title} {\bibinfo {title} {{Transfer learning for
  nonlinear dynamics and its application to fluid turbulence}},\ }\href@noop {}
  {\bibfield  {journal} {\bibinfo  {journal} {Physical Review E}\ }\textbf
  {\bibinfo {volume} {102}},\ \bibinfo {pages} {43301} (\bibinfo {year}
  {2020})}\BibitemShut {NoStop}%
\bibitem [{\citenamefont {Weiss}\ \emph {et~al.}(2016)\citenamefont {Weiss},
  \citenamefont {Khoshgoftaar},\ and\ \citenamefont {Wang}}]{weiss2016survey}%
  \BibitemOpen
  \bibfield  {author} {\bibinfo {author} {\bibfnamefont {K.}~\bibnamefont
  {Weiss}}, \bibinfo {author} {\bibfnamefont {T.~M.}\ \bibnamefont
  {Khoshgoftaar}},\ and\ \bibinfo {author} {\bibfnamefont {D.}~\bibnamefont
  {Wang}},\ }\bibfield  {title} {\bibinfo {title} {{A survey of transfer
  learning}},\ }\href@noop {} {\bibfield  {journal} {\bibinfo  {journal}
  {Journal of Big data}\ }\textbf {\bibinfo {volume} {3}},\ \bibinfo {pages}
  {1} (\bibinfo {year} {2016})}\BibitemShut {NoStop}%
\bibitem [{\citenamefont {do~Carmo}(1992)}]{do:carmo:1993}%
  \BibitemOpen
  \bibfield  {author} {\bibinfo {author} {\bibfnamefont {M.~P.}\ \bibnamefont
  {do~Carmo}},\ }\href@noop {} {\emph {\bibinfo {title} {{Riemannian
  Geometry}}}}\ (\bibinfo  {publisher} {Birkh{\"{a}}user Boston},\ \bibinfo
  {year} {1992})\BibitemShut {NoStop}%
\bibitem [{\citenamefont {Grigoryeva}\ and\ \citenamefont
  {Ortega}(2019)}]{RC9}%
  \BibitemOpen
  \bibfield  {author} {\bibinfo {author} {\bibfnamefont {L.}~\bibnamefont
  {Grigoryeva}}\ and\ \bibinfo {author} {\bibfnamefont {J.-P.}\ \bibnamefont
  {Ortega}},\ }\bibfield  {title} {\bibinfo {title} {{Differentiable reservoir
  computing}},\ }\href@noop {} {\bibfield  {journal} {\bibinfo  {journal}
  {Journal of Machine Learning Research}\ }\textbf {\bibinfo {volume} {20}},\
  \bibinfo {pages} {1} (\bibinfo {year} {2019})}\BibitemShut {NoStop}%
\bibitem [{\citenamefont {Jaeger}(2010)}]{jaeger2001}%
  \BibitemOpen
  \bibfield  {author} {\bibinfo {author} {\bibfnamefont {H.}~\bibnamefont
  {Jaeger}},\ }\href@noop {} {\emph {\bibinfo {title} {German National Research
  Center for Information Technology}}},\ \bibinfo {type} {Tech. Rep.}\
  (\bibinfo  {institution} {German National Research Center for Information
  Technology},\ \bibinfo {year} {2010})\BibitemShut {NoStop}%
\bibitem [{\citenamefont {Manjunath}\ and\ \citenamefont
  {Jaeger}(2013)}]{Manjunath:Jaeger}%
  \BibitemOpen
  \bibfield  {author} {\bibinfo {author} {\bibfnamefont {G.}~\bibnamefont
  {Manjunath}}\ and\ \bibinfo {author} {\bibfnamefont {H.}~\bibnamefont
  {Jaeger}},\ }\bibfield  {title} {\bibinfo {title} {{Echo state property
  linked to an input: exploring a fundamental characteristic of recurrent
  neural networks}},\ }\href {https://doi.org/10.1162/NECO_a_00411} {\bibfield
  {journal} {\bibinfo  {journal} {Neural Computation}\ }\textbf {\bibinfo
  {volume} {25}},\ \bibinfo {pages} {671} (\bibinfo {year} {2013})},\ \Eprint
  {https://arxiv.org/abs/1309.2848v1} {arXiv:1309.2848v1} \BibitemShut
  {NoStop}%
\bibitem [{\citenamefont {Manjunath}(2020)}]{manjunath:prsl}%
  \BibitemOpen
  \bibfield  {author} {\bibinfo {author} {\bibfnamefont {G.}~\bibnamefont
  {Manjunath}},\ }\bibfield  {title} {\bibinfo {title} {{Stability and
  memory-loss go hand-in-hand: three results in dynamics {\&} computation}},\
  }\href {https://doi.org/10.1098/rspa.2020.0563} {\bibfield  {journal}
  {\bibinfo  {journal} {To appear in Proceedings of the Royal Society London
  Ser. A Math. Phys. Eng. Sci.}\ ,\ \bibinfo {pages} {1}} (\bibinfo {year}
  {2020})},\ \Eprint {https://arxiv.org/abs/2001.00766} {arXiv:2001.00766}
  \BibitemShut {NoStop}%
\bibitem [{\citenamefont {Grigoryeva}\ and\ \citenamefont
  {Ortega}(2018{\natexlab{b}})}]{RC6}%
  \BibitemOpen
  \bibfield  {author} {\bibinfo {author} {\bibfnamefont {L.}~\bibnamefont
  {Grigoryeva}}\ and\ \bibinfo {author} {\bibfnamefont {J.-P.}\ \bibnamefont
  {Ortega}},\ }\bibfield  {title} {\bibinfo {title} {{Universal discrete-time
  reservoir computers with stochastic inputs and linear readouts using
  non-homogeneous state-affine systems}},\ }\href
  {http://arxiv.org/abs/1712.00754} {\bibfield  {journal} {\bibinfo  {journal}
  {Journal of Machine Learning Research}\ }\textbf {\bibinfo {volume} {19}},\
  \bibinfo {pages} {1} (\bibinfo {year} {2018}{\natexlab{b}})},\ \Eprint
  {https://arxiv.org/abs/1712.00754} {arXiv:1712.00754} \BibitemShut {NoStop}%
\bibitem [{\citenamefont {Lu}\ and\ \citenamefont
  {Bassett}(2020)}]{lu:bassett:2020}%
  \BibitemOpen
  \bibfield  {author} {\bibinfo {author} {\bibfnamefont {Z.}~\bibnamefont
  {Lu}}\ and\ \bibinfo {author} {\bibfnamefont {D.~S.}\ \bibnamefont
  {Bassett}},\ }\bibfield  {title} {\bibinfo {title} {{Invertible generalized
  synchronization: A putative mechanism for implicit learning in neural
  systems}},\ }\bibfield  {journal} {\bibinfo  {journal} {Chaos}\ }\textbf
  {\bibinfo {volume} {30}},\ \href {https://doi.org/10.1063/5.0004344}
  {10.1063/5.0004344} (\bibinfo {year} {2020})\BibitemShut {NoStop}%
\bibitem [{\citenamefont {Verzelli}\ \emph {et~al.}(2020)\citenamefont
  {Verzelli}, \citenamefont {Alippi},\ and\ \citenamefont
  {Livi}}]{Verzelli2020b}%
  \BibitemOpen
  \bibfield  {author} {\bibinfo {author} {\bibfnamefont {P.}~\bibnamefont
  {Verzelli}}, \bibinfo {author} {\bibfnamefont {C.}~\bibnamefont {Alippi}},\
  and\ \bibinfo {author} {\bibfnamefont {L.}~\bibnamefont {Livi}},\ }\bibfield
  {title} {\bibinfo {title} {{Learn to Synchronize, Synchronize to Learn}},\
  }\href {http://arxiv.org/abs/2010.02860} {\  (\bibinfo {year} {2020})},\
  \Eprint {https://arxiv.org/abs/2010.02860} {arXiv:2010.02860} \BibitemShut
  {NoStop}%
\bibitem [{\citenamefont {Hart}\ \emph {et~al.}(2021)\citenamefont {Hart},
  \citenamefont {Hook},\ and\ \citenamefont {Dawes}}]{allen:tikhonov}%
  \BibitemOpen
  \bibfield  {author} {\bibinfo {author} {\bibfnamefont {A.~G.}\ \bibnamefont
  {Hart}}, \bibinfo {author} {\bibfnamefont {J.~L.}\ \bibnamefont {Hook}},\
  and\ \bibinfo {author} {\bibfnamefont {J.~H.~P.}\ \bibnamefont {Dawes}},\
  }\bibfield  {title} {\bibinfo {title} {{Echo State Networks trained by
  Tikhonov least squares are L2($\mu$) approximators of ergodic dynamical
  systems}},\ }\href@noop {} {\bibfield  {journal} {\bibinfo  {journal}
  {Physica D: Nonlinear Phenomena}\ ,\ \bibinfo {pages} {132882}} (\bibinfo
  {year} {2021})}\BibitemShut {NoStop}%
\bibitem [{\citenamefont {Kocarev}\ and\ \citenamefont
  {Parlitz}(1996)}]{kocarev1996generalized}%
  \BibitemOpen
  \bibfield  {author} {\bibinfo {author} {\bibfnamefont {L.}~\bibnamefont
  {Kocarev}}\ and\ \bibinfo {author} {\bibfnamefont {U.}~\bibnamefont
  {Parlitz}},\ }\bibfield  {title} {\bibinfo {title} {{Generalized
  synchronization, predictability, and equivalence of unidirectionally coupled
  dynamical systems}},\ }\href@noop {} {\bibfield  {journal} {\bibinfo
  {journal} {Physical Review Letters}\ }\textbf {\bibinfo {volume} {76}},\
  \bibinfo {pages} {1816} (\bibinfo {year} {1996})}\BibitemShut {NoStop}%
\bibitem [{\citenamefont {Ceni}\ \emph
  {et~al.}(2020{\natexlab{a}})\citenamefont {Ceni}, \citenamefont {Ashwin},
  \citenamefont {Livi},\ and\ \citenamefont
  {Postlethwaite}}]{livi:multistability}%
  \BibitemOpen
  \bibfield  {author} {\bibinfo {author} {\bibfnamefont {A.}~\bibnamefont
  {Ceni}}, \bibinfo {author} {\bibfnamefont {P.}~\bibnamefont {Ashwin}},
  \bibinfo {author} {\bibfnamefont {L.}~\bibnamefont {Livi}},\ and\ \bibinfo
  {author} {\bibfnamefont {C.}~\bibnamefont {Postlethwaite}},\ }\bibfield
  {title} {\bibinfo {title} {{The echo index and multistability in input-driven
  recurrent neural networks}},\ }\href
  {https://doi.org/10.1016/j.physd.2020.132609} {\bibfield  {journal} {\bibinfo
   {journal} {Physica D: Nonlinear Phenomena}\ }\textbf {\bibinfo {volume}
  {412}},\ \bibinfo {pages} {132609} (\bibinfo {year} {2020}{\natexlab{a}})},\
  \Eprint {https://arxiv.org/abs/2001.07694} {arXiv:2001.07694} \BibitemShut
  {NoStop}%
\bibitem [{\citenamefont {Ceni}\ \emph
  {et~al.}(2020{\natexlab{b}})\citenamefont {Ceni}, \citenamefont {Ashwin},\
  and\ \citenamefont {Livi}}]{ceni:ashwin:paper1}%
  \BibitemOpen
  \bibfield  {author} {\bibinfo {author} {\bibfnamefont {A.}~\bibnamefont
  {Ceni}}, \bibinfo {author} {\bibfnamefont {P.}~\bibnamefont {Ashwin}},\ and\
  \bibinfo {author} {\bibfnamefont {L.}~\bibnamefont {Livi}},\ }\bibfield
  {title} {\bibinfo {title} {{Interpreting recurrent neural networks behaviour
  via excitable network attractors}},\ }\href
  {https://doi.org/10.1007/s12559-019-09634-2} {\bibfield  {journal} {\bibinfo
  {journal} {Cognitive Computation}\ }\textbf {\bibinfo {volume} {12}},\
  \bibinfo {pages} {330} (\bibinfo {year} {2020}{\natexlab{b}})},\ \Eprint
  {https://arxiv.org/abs/1807.10478} {arXiv:1807.10478} \BibitemShut {NoStop}%
\bibitem [{\citenamefont {Grigoryeva}\ \emph {et~al.}(2020)\citenamefont
  {Grigoryeva}, \citenamefont {Hart},\ and\ \citenamefont {Ortega}}]{RC19}%
  \BibitemOpen
  \bibfield  {author} {\bibinfo {author} {\bibfnamefont {L.}~\bibnamefont
  {Grigoryeva}}, \bibinfo {author} {\bibfnamefont {A.~G.}\ \bibnamefont
  {Hart}},\ and\ \bibinfo {author} {\bibfnamefont {J.-P.}\ \bibnamefont
  {Ortega}},\ }\bibfield  {title} {\bibinfo {title} {{Embedding chaos on
  manifolds with state-space systems}},\ }\href@noop {} {\bibfield  {journal}
  {\bibinfo  {journal} {In preparation}\ } (\bibinfo {year}
  {2020})}\BibitemShut {NoStop}%
\bibitem [{\citenamefont {Lorenz}(1963)}]{lorenz1963deterministic}%
  \BibitemOpen
  \bibfield  {author} {\bibinfo {author} {\bibfnamefont {E.~N.}\ \bibnamefont
  {Lorenz}},\ }\href@noop {} {\bibinfo {title} {{Deterministic nonperiodic
  flow}}} (\bibinfo {year} {1963})\BibitemShut {NoStop}%
\bibitem [{\citenamefont {Gonon}\ \emph {et~al.}(2020)\citenamefont {Gonon},
  \citenamefont {Grigoryeva},\ and\ \citenamefont {Ortega}}]{RC10}%
  \BibitemOpen
  \bibfield  {author} {\bibinfo {author} {\bibfnamefont {L.}~\bibnamefont
  {Gonon}}, \bibinfo {author} {\bibfnamefont {L.}~\bibnamefont {Grigoryeva}},\
  and\ \bibinfo {author} {\bibfnamefont {J.-P.}\ \bibnamefont {Ortega}},\
  }\bibfield  {title} {\bibinfo {title} {{Risk bounds for reservoir
  computing}},\ }\href@noop {} {\bibfield  {journal} {\bibinfo  {journal}
  {Journal of Machine Learning Research}\ }\textbf {\bibinfo {volume} {21}},\
  \bibinfo {pages} {1} (\bibinfo {year} {2020})}\BibitemShut {NoStop}%
\bibitem [{\citenamefont {Munkres}(2014)}]{Munkres:topology}%
  \BibitemOpen
  \bibfield  {author} {\bibinfo {author} {\bibfnamefont {J.}~\bibnamefont
  {Munkres}},\ }\href@noop {} {\emph {\bibinfo {title} {{Topology}}}},\
  \bibinfo {edition} {2nd}\ ed.\ (\bibinfo  {publisher} {Pearson},\ \bibinfo
  {year} {2014})\ p.\ \bibinfo {pages} {503}\BibitemShut {NoStop}%
\bibitem [{\citenamefont {Shapiro}(2016)}]{Shapiro:Farrago}%
  \BibitemOpen
  \bibfield  {author} {\bibinfo {author} {\bibfnamefont {J.~H.}\ \bibnamefont
  {Shapiro}},\ }\href@noop {} {\emph {\bibinfo {title} {{A Fixed-Point
  Farrago}}}}\ (\bibinfo  {publisher} {Springer International Publishing
  Switzerland},\ \bibinfo {year} {2016})\ p.\ \bibinfo {pages}
  {221}\BibitemShut {NoStop}%
\bibitem [{\citenamefont {Abraham}\ and\ \citenamefont
  {Robbin}(1967)}]{abraham:robbin}%
  \BibitemOpen
  \bibfield  {author} {\bibinfo {author} {\bibfnamefont {R.}~\bibnamefont
  {Abraham}}\ and\ \bibinfo {author} {\bibfnamefont {J.}~\bibnamefont
  {Robbin}},\ }\href@noop {} {\emph {\bibinfo {title} {{Transversal Mappings
  and Flows}}}}\ (\bibinfo  {publisher} {W. A. Benjamin},\ \bibinfo {address}
  {Inc},\ \bibinfo {year} {1967})\BibitemShut {NoStop}%
\bibitem [{\citenamefont {Hirsch}(1976)}]{Hirsch:book}%
  \BibitemOpen
  \bibfield  {author} {\bibinfo {author} {\bibfnamefont {M.~W.}\ \bibnamefont
  {Hirsch}},\ }\href@noop {} {\emph {\bibinfo {title} {{Differential
  Topology}}}}\ (\bibinfo  {publisher} {Springer Verlag},\ \bibinfo {year}
  {1976})\BibitemShut {NoStop}%
\end{thebibliography}


%

\end{document}